\documentclass[12pt]{amsart}

\usepackage{graphicx} 
\usepackage{amssymb,amsmath}
\usepackage{mathabx}
\usepackage[dvipsnames]{xcolor}
\usepackage{amsthm}
\usepackage{enumitem, csquotes}
\usepackage{fullpage}
\usepackage{hyperref}
\usepackage{cancel,tikz}
\usepackage{subcaption}

\DeclareMathOperator{\sgn}{sgn}

\DeclareMathOperator{\Rea}{Re}

\DeclareMathOperator{\supp}{supp}

\newtheorem{thm}{Theorem}[section]
\newtheorem{lemma}{Lemma}[section]
\newtheorem{remark}{Remark}[section]
\newtheorem{defn}{Definition}[section]
\newtheorem{prop}{Proposition}[section]
\newtheorem{cor}{Corollary}[section]

\newcommand{\p}{\partial}
\newcommand{\px}{\partial_x}
\newcommand{\pt}{\partial_t}
\newcommand{\pv}{\partial_v}
\newcommand{\la}{\langle}
\newcommand{\ra}{\rangle}

\newcommand{\ep}{\epsilon}

\newcommand{\R}{\mathbb{R}}

\newcommand{\wh}{\widehat}

\newcommand{\lp}{\left(}
\newcommand{\rp}{\right)}
\newcommand{\cd}{\cdot}

\newcommand{\lb}{\left[}
\newcommand{\rb}{\right]}

\newcommand{\bu}{\bar{u}}
\newcommand{\bv}{\bar{v}}

\newcommand{\bz}{\bar{z}}
\newcommand{\bg}{\bar{g}}

\newcommand{\up}{u^+}
\newcommand{\um}{u^-}

\renewcommand{\wp}{w^+}
\newcommand{\wm}{w^-}
\newcommand{\fp}{f^+}
\newcommand{\fm}{f^-}
\newcommand{\gp}{\gamma^+}
\newcommand{\gm}{\gamma^-}
\newcommand{\g}{\gamma}

\numberwithin{equation}{section}

\title{Dispersive decay for the modified Benjamin-Ono equation with small data}
\author{Allison Byars}
\date{}

\begin{document}
\maketitle

\begin{abstract}
    We study the modified Benjamin-Ono equation for small, localized initial data where no solitons arise, and prove a dispersive estimate.  Unlike the standard Benjamin-Ono equation, the modified version is not known to be completely integrable, so we cannot rely on integrable systems tools. First, to manage the nonlocal Hilbert transform in the linear part of the equation, we project the equation to  positive and negative frequencies. We then combine wave packet testing of Ifrim–Tataru, the Klainerman–Sobolev vector field method, and a bootstrap argument. 
\end{abstract}

\tableofcontents

\section{Introduction\label{Intro mbo section}}
We consider the Cauchy problem for the modified Benjamin-Ono equation (mBO), which is given by
\begin{equation} \label{mBO} \tag{mBO}
\begin{cases}
     u_t+H\partial_x^2u =\partial_{x}(|u|^2u)\\
     u(0,x)=u_0(x),
\end{cases}
\end{equation}
where $u: \R\times \R \to \mathbb{C}.$  We note that usually for this equation, only real valued solutions are considered.  Everything in this paper also holds for the real valued case. Indeed, \eqref{mBO} is invariant under complex conjugation, since $\overline{H\px^2 u}=H\px^2\bu$ and $\overline{|u|^2u}=|\bu|^2\bu$. In particular real data lead to real solutions, where $|u|^2u=u^3$. 

The modified Benjamin-Ono equation is sometimes written as
\begin{equation*}
    u_t+H\partial_x^2u =\pm \partial_{x}(|u|^2u).
\end{equation*}
and is also sometimes called cubic Benjamin-Ono. We fix the $+$ sign throughout.
$H$ is the Hilbert transform, given by
\begin{align}\label{Hilbert transform defn}
    Hu(x)=\frac{1}{\pi} p.v. \int_{-\infty}^\infty \frac{u(y)}{x-y}\, dy,
\end{align}
which has symbol $-i\sgn \xi$. Therefore, the dispersion relation for \eqref{mBO} is $a(\xi)=\xi|\xi|.$
\eqref{mBO} admits the following scaling symmetry: If $u$ is a solution, then for $\lambda>0$ so is
\begin{align*}
    u_\lambda(t,x)=\lambda^{\frac{1}{2}}u( \lambda^2 t,\lambda x).
\end{align*}
From this scaling, we can find the critical Sobolev space, which is the space $\dot{H^s}$ where
\[ \|u_\lambda\|_{\dot{H^s}}=\|u\|_{\dot{H^s}}.\]
In this case, the critical Sobolev space is $L^2.$ This is the natural threshold to aim for when proving well-posedness. From the scaling, we can also get the scaling vector field
\begin{align}\label{L defn mbo}
    L=x -2tH\px.
\end{align}
This commutes with the linear part of the equation since
\begin{align*}
     [\pt+H\px^2, L]u=&[H,x]\px^2 u=-\frac{1}{\pi}\int_\R \px^2 u(y) \, dy =0.
\end{align*}
This commutator can be seen in Lemma \ref{Hilbert Transform Lemma}.
It also has the property that
\begin{align*}
    [\px, L]=&1,
\end{align*}
and it is the pushforward of $x$ along the linear flow.  In other words, 
\[L(t)=e^{-tH\px^2}\,x\,e^{tH\px^2}.\]
It is not known whether  \eqref{mBO} is completely integrable.  For example, there is no known Lax pair or inverse scattering transform.  It is also not known whether it has infinitely many conserved quantities. The three known conserved quantities are the integral of the solution:
\begin{align*}
    \frac{d}{dt}\int u \, dx =0,
\end{align*}
the $L^2-$ mass:
\begin{align*}
    \frac{d}{dt}\int |u|^2 \, dx =0,
\end{align*}
and the Hamiltonian:
\begin{align*}
    \frac{d}{dt}\int \lp \bu Hu_x- \frac{1}{2}|u|^4 \rp\, dx =0.
\end{align*}
To verify that this is the Hamiltonian, we let $\mathcal{H}(u)=-\frac{1}{2}\int \lp \bu Hu_x- \frac{1}{2}|u|^4 \rp\, dx.$  Then,
\begin{align*}
 \delta \mathcal{H}(u)=-Hu_x+u^2 \bu,
\end{align*}
and for $\mathcal{J}=\px,$ we have $u_t=\mathcal{J}(\delta \mathcal{H}(u)).$  Here, $\delta\mathcal H$ denotes the functional (variational) derivative with respect to the real inner product $\Rea\int f\bg.$

\textbf{History of \eqref{mBO}:}

The Benjamin-Ono family of equations models internal waves at the interface of two deep fluids, where the Hilbert transform captures nonlocal dispersive effects. The standard Benjamin–Ono (BO) equation was derived independently by Benjamin \cite{benjamin1967internal} and Ono \cite{ono1975algebraic}, and accounts for quadratic wave interactions. The modified version \eqref{mBO} shares the same linear dispersive structure as the classical BO equation, but it replaces the quadratic nonlinearity with a cubic one to model stronger wave interactions.

The BO equation is a completely integrable model, meaning it has an inverse scattering transform, a Lax pair, and an infinite number of conservation laws.  However, its modified version mBO is not known to be integrable. The potential loss of integrability introduced by the cubic term makes the low-regularity theory for mBO substantially more challenging.

We now discuss the history of the well-posedness for mBO in $H^s$.  First, it is $L^2$ critical, i.e. critical at $s=0.$ 
For the following results, only real valued solutions were considered. 
In 2003, Kenig and Koenig proved local well-posedness for $s\geq 1$ \cite{kenig2003local}. A year later, Molinet and Ribaud lowered this to $s>\frac{1}{2}$, using a gauge transformation \cite{molinet2004well}. 
In 2006, Kenig and Takaoka proved the global well-posedness in the case $s=\frac{1}{2},$ again relying on a gauge transformation, but this time combined with Littlewood-Paley \cite{kenig2006global}.  They also proved that this threshold is sharp, in the sense that the data-to-solution map fails to be uniformly continuous below $s=\frac{1}{2}$.

While uniform continuity fails below this point, one can still go for a weaker notion of well-posedness, requiring only continuity of the solution map. From this viewpoint, the entire range $0 \leq s < \frac{1}{2}$ remains open. There has been some progress in this direction. In 2008, Guo studied the problem with complex data, and proved  local well-posedness in the case $s\geq \frac{1}{2}$ without using a gauge transformation \cite{guo2008local}.  In the same paper, he proved global well-posedness for real data in $H^s$ for $s\geq \frac{1}{2}$, and a priori bounds for $s>1/4.$ There are also numerous works on more generalized versions of the equation  \cite{ai2024dispersion,guo2012local,kenig1994generalized}.

The long-time behavior of \eqref{mBO} has been studied less. When working towards proving the Soliton Resolution conjecture, it is natural to first consider the case where there are no solitons. In 1999, Hayashi and Naumkin \cite{hayashi1999large} proved a dispersive estimate under the assumption that the initial data is small in $H^{3,0} \cap H^{1,2}$, where the  $H^{s,m}$ norms are defined as 
\begin{align*}
    \|f\|_{H^{s,m}}:=\|(1+|x|^2)^{\frac{m}{2}}(1-\px^2)^{\frac{s}{2}}f\|_{L^2_x}.
\end{align*}
In 2001, they improved on this slightly by lowering the assumption to initial data in $H^{2,0}\cap H^{1,1}$ \cite{hayashi2001reduction}. The goal of this paper will be to lower these assumptions even more, with only $u_0\in H^1$ and $xu_0\in H^1$, being required. This is equivalent to $u_0\in H^{1,1}.$ We also generalize the result to the case where $u$ is complex valued. More specifically, the main theorem is as follows.
\begin{thm}
\label{mbo thm}
    Let $u$ be a solution to \eqref{mBO} with small and localized data, i.e. for $\epsilon \ll 1,$
    \begin{equation}
        \|xu_0\|_{H^1}+\|u_0\|_{H^1}\leq \epsilon.
        \label{initials mbo}
    \end{equation}
Then we have the following bound on $Lu$,
\begin{equation}
    \|Lu\|_{H^1_x}\lesssim \epsilon\, \langle t \rangle^{C\epsilon^2},
    \label{Lubound1 mbo}
\end{equation}
and the dispersive bounds
\begin{equation*}
\|u(t,x)\|_{L^\infty_x}+\|u_x(t,x)\|_{L^\infty_x} \lesssim \epsilon \, \langle t\rangle^{-\frac{1}{2}} \text{ for all } t\in \R.
\end{equation*}
\end{thm}
Note: In the proof, we show both bounds for $t\geq 1.$  The bounds hold in the range $0\leq t\leq 1$  by a short bootstrap argument, which will be given in section \ref{Close bootstrap section}.  We then obtain the case $t<0$ by the symmetry $(t,x)\mapsto(-t,-x)$.

\subsection{Solitons} 

We now briefly study the solitons for \eqref{mBO} to show that Theorem \ref{mbo thm} is considering the case of no solitons.  The solitons are given in \cite{roudenko2021dynamics}, but we derive them here for completeness.

Let $q(t,x)=c^a Q(c(x + ct)) $ be a soliton. Then, in order for $q$ to satisfy \eqref{mBO},
we would need
\begin{align*}
&\partial_t q + H \partial_x^2 q - \partial_x (|q|^2 q)
=\partial_x \left[ c^{a+1} Q(cx + c^2 t) + c^{a+1} HQ'(cx + c^2 t) -c^{3a} |Q|^2 Q \right]=0.
\end{align*} 
In order to factor out the constants, we take $a=\frac{1}{2}.$ Then,
as long as $Q$ satisfies 
\begin{equation}\label{Q eqn}
    Q+HQ'-|Q|^2Q=0,
\end{equation}
we have
\begin{align*}
    \partial_t q + H \partial_x^2 q - \partial_x (|q|^2 q)=0.
\end{align*}
As in \cite{boyd2012numerical,roudenko2021dynamics} let $Q$ be the unique positive, even, decreasing for $x>0,$ vanishing at infinity solution to \eqref{Q eqn}. $Q$ has the property that $Q(x) = O(|x|^{-2})$.  The family of solitons is then given by
\[q(t,x) = c^{\frac{1}{2}} Q(c(x + ct)), \text{ for } c> 0.\]
\bigskip
Next, we compare the norms of $q_0(x) = c^{\frac{1}{2}} Q(cx)$, the ground-state soliton profile with wave speed $c > 0$, where $Q(y)$ satisfies $Q(y) \sim O(|y|^{-2})$ and $Q'(y) \sim O(|y|^{-3})$ as $|y| \to \infty$. First, we compute the $L^2$ norm of the derivative profile $q_0'(x) = c^{\frac{3}{2}} Q'(cx)$:
\begin{align*}
\|q_0'\|_{L^2}^2 &= \int_{-\infty}^{\infty} |c^{\frac{3}{2}} Q'(cx)|^2 \, dx \\
&= c^2 \int_{-\infty}^{\infty} |Q'(y)|^2 \, dy = c^2 \|Q'\|_{L^2}^2.
\end{align*}
The localized $L^2$ norm is:
\begin{align*}
\|x q_0\|_{L^2}^2 &= \int_{-\infty}^{\infty} x^2 |c^{\frac{1}{2}} Q(cx)|^2 \, dx \\
&= \frac{1}{c^2} \|y Q\|_{L^2}^2=O(c^{-2}).
\end{align*}
Then
\[\|q_0'\|_{L^2}=O(c), \quad \|x q_0\|_{L^2}=O(c^{-1}),\]
and so
\[
\|x q_0 \|_{L^2} \cdot \|q_0' \|_{L^2} \gtrsim 1. \quad 
\]
Therefore, $q$ cannot satisfy both $\|x q_0 \|_{L^2}\leq \epsilon$ and $\|q_0' \|_{L^2} \leq \epsilon$ for $\epsilon$ small enough.

The solitons above are real valued. Since \eqref{mBO} is invariant under the phase rotation $u\mapsto e^{i\theta}u$, for complex solutions one also has the solitons $e^{i\theta}q$, for $\theta\in\R$. These have the same norms, so the conclusion is unchanged.

\subsection{Acknowledgments}
The author would like to thank Mihaela Ifrim for proposing the problem and for many helpful discussions.  The author was supported in part by NSF grant DMS-2037851 in Spring 2026, NSF DMS-2348908 and a Vilas Associate Fellowship in Summer 2026, and NSF RTG DMS-2135998 for Fall 2026.

\section{Preliminaries}
In this section, we give some basic facts and definitions that will be used throughout the paper.
First, we are using the following definition of the Fourier transform:
\begin{align*}
    \hat{f}(\xi):=\frac{1}{\sqrt{2\pi}}\int_{\R}e^{-ix\xi}f(x)\, dx.
\end{align*}
Next, because of the Hilbert transform in the equation, we will often need to project to positive and negative  frequencies, and consider each case separately.  

\begin{defn} \label{LP defn}
Define the positive and negative frequency projections $P^+, P^-$ as follows
    \begin{align*}
        \wh{P^+f}(\xi)=&\chi_{\xi>0}(\xi)\wh{f}(\xi)\\
     \wh{P^-f}(\xi)=&\chi_{\xi<0}(\xi)\wh{f}(\xi),
    \end{align*}
    Notation: We will sometimes write $f^+=P^+f$ and $f^-=P^-f.$
    Note: $\widebar{f^+}=P^-\widebar{f}$ and $\widebar{f^-}=\widebar{f}^+$ 
since, using the identity $\widehat{\overline f}(\xi)=\overline{\widehat f(-\xi)}$, we have
\begin{align*}
    \wh{\overline{f^+}}(\xi)=\overline{\wh{\fp}(-\xi)}=\chi_{\{\xi>0\}}(-\xi)\overline{\wh{f}(-\xi)}=\chi_{\{\xi<0\}}(\xi)\overline{\wh{f}(-\xi)}=\wh{(\overline{f})^-}(\xi).
\end{align*}
\end{defn}

\begin{thm}[Plancherel Theorem]\label{plancherel}
    Let $f,g$ be two $L^2$ functions. Then,
\begin{align*}
    \int_{-\infty}^\infty f(x) \widebar{g(x)}\, dx=   \int_{-\infty}^\infty \wh{f}(\xi) \overline{\wh{g}(\xi)}\, d\xi.
\end{align*}
This also implies
\begin{align*}
    \|f\|_{L^2_x}=\|\wh{f}\|_{L^2_\xi}.
\end{align*}
\end{thm}

\begin{lemma}
    $P^+$ is self adjoint. In other words, for two functions $f,g$, $\la P^+f,g \ra_{L^2}=\la f,P^+ g \ra_{L^2}$.
\end{lemma}
\begin{proof}
Using Plancherel, we see
    \begin{align*}
        \la P^+f,g \ra_{L^2}=&\int  P^+f(x)\widebar{g(x)} \, dx=\int \wh{P^+f}(\xi) \overline{\wh{g}(\xi)}\, d\xi\\
        =&\int \chi_{+}(\xi) \wh{f}(\xi) \overline{\wh{g}(\xi)}\, d\xi\\
        =&\int  \wh{f}(\xi) \overline{\chi_{+}(\xi)\wh{g}(\xi)}\, d\xi\\
        =&\int  \wh{f}(\xi) \overline{\wh{P^+ g}(\xi)}\, d\xi\\
        =&\int f(x) \overline{P^+g(x)}\, dx=     \la f,P^+g \ra_{L^2}.
    \end{align*}
\end{proof}

\begin{lemma}[Hilbert Transform Properties] \label{Hilbert Transform Lemma}
    For $H$ as defined in \eqref{Hilbert transform defn}, we have the following properties: 
  \begin{itemize}
    \item $H^{-1}=-H$
    \item $H^2=-1$
    \item $\p^k H(u)=H(\p^k u)$
    \item $\widebar{Hu}=H \widebar{u}$
    \item $\wh{Hu}(\xi)=-i\sgn \xi \cdot \wh{u}(\xi)$
    \item $[H,x]u=-\frac{1}{\pi}\int u(y)\, dy$
\end{itemize}
Also, for any two functions $u,v\in L^2,$ we have
\begin{itemize}
    \item $\langle Hu,v\rangle_{L^2} =\langle u,-Hv\rangle_{L^2} $
    \item $\|Hu\|_{L^2}=\|u\|_{L^2}$.
\end{itemize}
Additionally, we can write $iH=P^+-P^-$ and we have $P^+H=-iP^+, P^-H=iP^-.$
We can also write $P^+, P^-$ in terms of the Hilbert transform as
\begin{align*}
    P^+:=&\frac{1}{2}(I+iH)\\
 P^-:=&\frac{1}{2}(I-iH).
\end{align*}
\end{lemma}

\begin{lemma}
  Let $L$ be the scaling vector field \eqref{L defn mbo}. Then  $L$ is self adjoint. In other words, for two Schwartz functions $f,g$, $\la Lf,g \ra_{L^2}=\la f,L g \ra_{L^2}$.
\end{lemma}

\begin{proof}
Using the properties of the Hilbert transform above, we have
\begin{align*}
    \la Lf, g\ra_{L^2}=&\int Lf \cdot \bg =\int (x-2tH\px)f \cdot \bg \\
    =& \int f \cdot x\bg +2t\px f \cdot H \bg\\
 =& \int f \cdot x\bg -2t f \cdot \px\overline{H g}\\
     =&\int f \widebar{Lg}= \la f, Lg\ra_{L^2}.
\end{align*}
\end{proof}

The price to pay for the sharp cutoff is that $P^\pm$ do not commute with multiplication by $x$.  Instead, the commutator is a constant.
\begin{lemma}[Commutators]\label{commutator lemma}
For $f\in L^2$ with $xf\in L^2$,
\begin{align}\label{x P commutator}
    [x,P^\pm]f=\pm\frac{i}{2\pi}\int_\R f\,dy.
\end{align}
Consequently, if $u\in C([0,T];H^1)$ solves \eqref{mBO} with $u_0\in H^{1,1}$, then for $Lu \in L^2$,
\begin{align}\label{L P commutator}
    Lu^\pm=P^\pm Lu\pm c,\qquad c:=\frac{i}{2\pi}\int_\R u_0\,dx,\qquad |c|\leq\frac{1}{2\pi}\|u_0\|_{L^1}\leq \|xu_0\|_{L^2}+\|u_0\|_{L^2} .
\end{align}
Moreover, 
\begin{align*}
   L(u^\pm_x)=P^\pm L(u_x).
\end{align*}

\end{lemma}
\begin{proof}
For the first commutator, we take the  Fourier transform, 
\begin{align*}
    (xP^+f)^\wedge(\xi)=i\p_\xi\lp\chi_{\{\xi>0\}}(\xi)\wh f(\xi)\rp=i\,\wh f(0)\,\delta_0+\chi_{\{\xi>0\}}(\xi)\,i\p_\xi\wh f=i\wh f(0)\delta_0+\lp P^+(xf)\rp^\wedge(\xi).
\end{align*}
Therefore
\begin{align*}
    [x,P^+]f=\frac{i}{2\pi}\int f.
\end{align*}
 The case of $P^-$ is identical with the opposite sign. For \eqref{L P commutator}, since $H, \px$ commute with $P^\pm,$ we have
 \begin{align*}
     [L,P^\pm]u=[x,P^\pm]u=\pm \frac{i}{2\pi}\int u \, dy =\pm \frac{i}{2\pi}\int u_0 \, dy,
 \end{align*}
since the integral of the solution is conserved.

The bound for $|c|$ follows from $\|u_0\|_{L^1}\leq\|\la x\ra^{-1}\|_{L^2}\|\la x\ra u_0\|_{L^2}\lesssim \|xu_0\|_{L^2}+\|u_0\|_{L^2}$. The last statement follows since
\[L(u^\pm_x)=P^\pm L(u_x)+\int u_x\,dx=P^\pm L(u_x).\]

\end{proof}
In particular, this Lemma tells us that $Lu^\pm\notin L^2$ (unless $c=0$), but $P^\pm Lu\in L^2$ with  $\|P^\pm Lu\|_{L^2}\leq\|Lu\|_{L^2}$. Moreover, if $Lu\in H^1$ then $\|Lu^\pm\|_{L^\infty}\leq\|Lu\|_{H^1}+|c|$.

The constant $c$ is a conserved quantity of size $O(\ep)$. It appears in the analysis only through the identity \eqref{L P commutator}, and it will be tracked explicitly.

Next, for some energy estimates, we will need a nonlinear version of Gr\"onwall's inequality. 
\begin{lemma}[Nonlinear Gr\"onwall]\label{NL Gronwall}
    For some functions $f,p,$ and $B,$ if 
    \begin{align*}
        \frac{d}{dt}f(t)^2\leq p(t) f(t)^2+ B(t)f(t),
    \end{align*}
    then
    \begin{align*}
          f(t)\leq& f(t_0) \cdot \exp\lp\int_{t_0}^t \frac{p(s)}{2}\, ds\rp+\frac{1}{2}\int_{t_0}^t B(s) \exp\lp\int_{s}^t \frac{p(r)}{2}\, dr\rp \, ds.
    \end{align*}
    In particular, if $p(t)=C_1t^{-1}$ and $B(t)=C_2t^{-1},$ where $C_1, C_2>0,$ then
    \begin{align*}
             |f(t)|\leq& |f(t_0)|\lp \frac{t}{t_0} \rp^{\frac{C_1}{2}} + \frac{C_2}{C_1}\left[\lp \frac{t}{t_0} \rp^{\frac{C_1}{2}}-1\right].
    \end{align*}
\end{lemma}
See \cite{byars2026modscatteringAC} for the proof.

\section{Linear Equation \label{linear sec}}
We first study the linear equation associated to \eqref{mBO}, in order to gain insights on the dynamics for the nonlinear version.  Here, the linear equation is the linear Benjamin-Ono equation
\begin{align}\label{linear eqn} \tag{Lin}
    u_t+H\px^2 u=0.
\end{align}
We can solve this by taking the Fourier Transform.  Let $t>0,$ then we have 
\begin{align*}
    \wh{u_t}(\xi)+&(-i\sgn(\xi))\cd (i\xi)^2 \wh{u}(\xi)=0\\
     {u}(t,x)=& u_0(x) * (e^{-i\sgn \xi \cd \xi^2t})^\vee(x).
\end{align*}
Then, since $\sgn \xi \cd \xi^2=|\xi|\xi$,
\begin{align*}
    (e^{-i\sgn \xi \cd \xi^2t})^\vee (x)=& \int e^{i x\cd \xi}e^{-i\sgn \xi \cd \xi^2t} \, d\xi\\
    =& t^{-\frac{1}{2}}\int e^{i x\cd t^{-\frac{1}{2}}\xi}e^{-i| \xi| \cd \xi} \, d\xi\\
    =&t^{-\frac{1}{2}}A(t^{-\frac{1}{2}}x),
\end{align*}
where
\begin{align*}
    A(x)=\int e^{i x\cd \xi}e^{-i| \xi| \cd \xi} \, d\xi.
\end{align*}
It remains to show that $A$ is uniformly bounded, which can be done by a stationary phase argument. 
\begin{align*}
        A(x)=&\int e^{i x\cd \xi}e^{-i| \xi| \cd \xi} \, d\xi\\
        =&\int_{-\infty}^0 e^{i x\cd \xi}e^{i \xi^2} \, d\xi+\int_0^\infty e^{i x\cd \xi}e^{-i \xi^2} \, d\xi\\
        =&\int_{-\infty}^0 e^{i (\xi+\frac{x}{2})^2}e^{-i \frac{x^2}{4}} \, d\xi  +   \int_0^\infty e^{-i(\xi-\frac{x}{2})^2}e^{i \frac{x^2}{4}} \, d\xi\\
          =&e^{-i \frac{x^2}{4}} \int_{-\infty}^0 e^{i (\xi+\frac{x}{2})^2}\, d\xi  +   e^{i \frac{x^2}{4}} \int_0^\infty e^{-i(\xi-\frac{x}{2})^2}\, d\xi\\
         =& e^{-i \frac{x^2}{4}}\int_{-\infty}^{\frac{x}{2}}e^{i\eta^2}\,d\eta+e^{i \frac{x^2}{4}}\int_{-\frac{x}{2}}^{\infty}e^{-i\eta^2}\,d\eta
\end{align*}
Then $\|A\|_\infty\lesssim1$, and we obtain
\begin{align*}
    u(t,x)=t^{-\frac{1}{2}} u_0(x)*A(t^{-\frac{1}{2}}x).
\end{align*}
By Young's inequality, this gives $\|u\|_{L^\infty}\lesssim t^{-\frac{1}{2}}\|A\|_{L^\infty}\|u_0\|_{L^1}.$ This gives the dispersive decay rate for the linear solution as $t^{-\frac{1}{2}}.$  Therefore, in the case of no solitons, the goal will be to prove this same decay rate for the nonlinear solution.

\subsection{\texorpdfstring{Linear Analysis using the vector field $L$}{Linear Analysis using the vector field L}}

Since $L$ commutes with the linear equation, we have, for $u$ a solution to \eqref{linear eqn},
\begin{align*}
    (\pt+H\px^2)Lu=0.
\end{align*}

The first idea might be to try and do the analysis as was done for NLS \cite{cite:NLS} and DNLS \cite{byars2024globaldynamicssmalldata}. However, the Hilbert transform poses a problem, which we will discover in the next section. Usually, we can obtain the dispersive estimate in another way by using the vector field. This is true for NLS and DNLS.  However, for \eqref{mBO}, the way we obtain this bound is different, and requires separation of the positive and negative frequency terms. First, we need a few definitions.

\begin{defn}
    Define $S(t) = e^{-t H \partial_x^2}$ to be the linear propagator for \eqref{mBO}.  Then we can write $(S(t) P^\pm f)(x) = (K^\pm (t,\cdot) * f) (x)$
    where the kernels are
\begin{align*}
K^+(t, x) =& \frac{1}{2\pi} \int_0^\infty e^{i x \xi - i t \xi^2} \, d\xi\\
    K^-(t, x) =& \frac{1}{2\pi} \int_{-\infty}^0 e^{i x \xi + i t \xi^2} \, d\xi.
\end{align*}
\end{defn}

\begin{lemma}
Let $S(t) = e^{-t H \partial_x^2}$ be as defined above.
The kernel of $S(t) P^\pm$ admits the scaling representation
\[
K^\pm(t, x) = t^{-\frac{1}{2}} A^\pm\left(t^{-\frac{1}{2}} x\right),
\]
where $A^\pm \in L^\infty(\mathbb{R})$ are partial Fresnel integrals. Consequently,
\[
\|K^\pm(t, \cdot)\|_{L^\infty} \lesssim t^{-\frac{1}{2}}.
\]
\end{lemma}

\begin{proof}
We first consider the positive frequency case. We have
\[
K^+(t, x) = \frac{1}{2\pi} \int_0^\infty e^{i x \xi - i t \xi^2} \, d\xi.
\]
Making the change of variables $\eta = t^{1/2} \xi$, $d\xi = t^{-1/2} d\eta$, we obtain
\[
K^+(t, x) = t^{-\frac{1}{2}} \frac{1}{2\pi} \int_0^\infty e^{i (t^{-1/2} x) \eta - i \eta^2} \, d\eta = t^{-\frac{1}{2}} A^+\left(t^{-\frac{1}{2}} x\right),
\]
where
\[
A^+(y) := \frac{1}{2\pi} \int_0^\infty e^{i y \eta - i \eta^2} \, d\eta.
\]
The function $A^+(y)$ is a partial Fresnel integral. Then by completing the square in the exponent and integrating by parts, $A^+$ is uniformly bounded on $\mathbb{R}$, so $A^+ \in L^\infty(\mathbb{R})$.

Similarly, for negative frequencies, substituting $\eta = -t^{1/2} \xi > 0$ yields
\[
K^-(t, x) = \frac{1}{2\pi} \int_{-\infty}^0 e^{i x \xi + i t \xi^2} \, d\xi = t^{-\frac{1}{2}} A^-\left(t^{-\frac{1}{2}} x\right), \quad A^-(y) = \overline{A^+(y)} \in L^\infty(\mathbb{R}).
\]
In each case we obtain the bound $\|K^\pm(t, \cdot)\|_{L^\infty} \le \|A^\pm\|_{L^\infty} t^{-1/2} \lesssim t^{-1/2}$.
\end{proof}

\begin{prop}\label{vector field estimate general f} Let $f\in L^2$ with $Lf\in L^2.$ Then $f^\pm$ satisfy
\begin{equation}\label{eq:disper_bound_pm}
\|f^\pm\|_{L^\infty_x} \lesssim t^{-\frac{1}{2}} \|f\|_{L^2}^{\frac{1}{2}} \|Lf\|_{L^2}^{\frac{1}{2}}.
\end{equation}
In particular, we apply this to $f=u, f=u_x.$
\end{prop}

\begin{proof}
Define the profile $g(t, x) := S(-t) f(t, x) = e^{t H \partial_x^2} f(t, x)$. Since $S(t)$ is unitary on $L^2(\mathbb{R})$, we have
\[
\|g\|_{L^2} = \|f\|_{L^2}.
\]
Since $L$ is the pushforward of $x$ along the linear flow, we have the identity $S(-t) L S(t) = x$, and so $x g = S(-t) L f$. This gives
\[
\|x g\|_{L^2} = \|S(-t) L f\|_{L^2} = \|Lf\|_{L^2}.
\]
Expressing $f^\pm = S(t) P^\pm g$ as a convolution with the kernel $K^\pm(t, x)$:
\[
f^\pm(x) = (K^\pm(t, \cdot) * g)(x) = \int_{\mathbb{R}} K^\pm(t, x - y) g(y) \, dy.
\]
Let $R>0$ be some positive constant which will be chosen later. We bound $f^\pm$ by Young's inequality and split up the integral:
\[
|f^\pm(x)| \le \|K^\pm(t, \cdot)\|_{L^\infty} \left( \int_{|y| \le R} |g(y)| \, dy + \int_{|y| > R} |g(y)| \, dy \right).
\]
We bound each integral by Cauchy--Schwarz:
\begin{align*}
     \int_{|y| \le R} |g(y)| \, dy \le& \left( \int_{|y| \le R} 1 \, dy \right)^{\frac{1}{2}} \|g\|_{L^2} = \sqrt{2R} \, \|g\|_{L^2}\\
      \int_{|y| > R} \frac{1}{|y|} |y g(y)| \, dy \le& \left( \int_{|y| > R} \frac{1}{y^2} \, dy \right)^{\frac{1}{2}} \|x g\|_{L^2} = \sqrt{\frac{2}{R}} \, \|x g\|_{L^2}.
\end{align*}
Combining these bounds with $\|K^\pm(t, \cdot)\|_{L^\infty} \lesssim t^{-1/2}$ yields
\[
|f^\pm(x)| \lesssim t^{-\frac{1}{2}} \left( \sqrt{R} \|g\|_{L^2} + \frac{1}{\sqrt{R}} \|x g\|_{L^2} \right).
\]
Setting $R = \frac{\|x g\|_{L^2}}{\|g\|_{L^2}}$ gives
\[
\|f^\pm\|_{L^\infty_x} \lesssim t^{-\frac{1}{2}} \|g\|_{L^2}^{\frac{1}{2}} \|x g\|_{L^2}^{\frac{1}{2}} = t^{-\frac{1}{2}} \|f\|_{L^2}^{\frac{1}{2}} \|Lf\|_{L^2}^{\frac{1}{2}}.
\]
\end{proof}

\begin{cor}
 For solutions $u$ to \eqref{linear eqn} with $\|u_0\|_{L^2}+\|xu_0\|_{L^2}\leq \epsilon$, we have the following dispersive bound:
    \begin{equation*}
    \|u(t,x)\|_{L^\infty_x}\lesssim \epsilon t^{-\frac{1}{2}}.
\end{equation*}
\end{cor}
\begin{proof}
    Apply the previous proposition to $\up, \um$ along with conservation of mass for $u$ and $Lu.$ 
\end{proof}

\section{Nonlinear Analysis \label{nonlinear sec}}
Although the vector field $L$ can be used to prove the dispersive estimate in the linear case, this does not translate well to the nonlinear case, because the vector field does not play nicely with the nonlinearity.  Therefore, as for DNLS in \cite{byars2024globaldynamicssmalldata}, we will use the testing by wave packets method of Ifrim and Tataru \cite{cite:NLS}. 

This will be significantly different from the DNLS case for a few reasons, all stemming from the Hilbert transform in the equation. First, when defining the wave packets, the dispersion relation must be twice differentiable.  For \eqref{mBO}, the dispersion relation is $a(\xi)=\xi |\xi|,$ which presents problems.  To solve this, we split up the equation by projecting onto positive, and then negative frequencies.  We then perform the wave packets analysis for the positive frequency projection of the solution, $P^+u=\up,$ and the negative projection $P^-u=\um$ separately.  

The next challenge comes when attempting to obtain energy estimates for $L u$.  This is again because of the Hilbert transform in $L.$ Now, $L$ does not have the same properties that it did for DNLS, so we cannot obtain the estimate directly.  The difficulties here and how to overcome them will be seen in more detail in section \ref{Vector field bnds section}.

The last few differences between the present analysis and the case of DNLS \cite{byars2024globaldynamicssmalldata} come from the asymptotic equation.  First, having the positive and negative frequency projections makes the remainder terms far more complicated to bound, and requires a careful analysis.  In addition, the asymptotic equation admits terms which are unbounded in the energy.  Therefore, we must perform a normal form analysis in order to deal with these nonresonant terms. 

The outline of the proof of Theorem \ref{mbo thm} is as follows.

\textbf{I. The Paradifferential equation}
As mentioned above, the Hilbert transform necessitates projecting the equation to positive and negative frequencies. We then do a wave packets analysis on each equation separately.  By applying $P^+$  to \eqref{mBO}, we get

\begin{align} \label{positive freq proj mbo}
    \pt u^+-i\px^2 u^+=P^+\px(|u|^2u).
\end{align}
Similarly, projecting to negative frequencies yields
\begin{align}\label{negative freq proj mbo}
    \pt u^-+i\px^2 u^-=P^-\px(|u|^2u).
\end{align}

\textbf{II. Define the wave packets}
We will perform the wave packets analysis by testing both $\up$ and $\um$ against their corresponding wave packets. Therefore, we define both the positive and negative frequency wave packets as follows.  
\begin{align*}
    \Phi_v^+=e^{i\frac{x^2}{4t}}\chi\lp \frac{x-tv}{t^{\frac{1}{2}}}\rp, \,\,\, \Phi_v^-(t,x)=e^{-i\frac{x^2}{4t}}\chi\lp \frac{x-tv}{t^{\frac{1}{2}}}\rp,
\end{align*}
where $\chi$ is a real, even bump function with $\int \chi(y)\, dy =1.$ The wave packets are approximate solutions to the positive (respectively negative) frequency linear equation:
\begin{align*}
    \pt \Phi_v^+-i\px^2\Phi_v^+
    =&\frac{e^{i\frac{x^2}{4t}}}{2t} \px \lb (x-tv)\chi\lp \frac{x-tv}{t^{\frac{1}{2}}}\rp -2it^{\frac{1}{2}} \chi'\lp \frac{x-tv}{t^{\frac{1}{2}}}\rp  \rb\\
    =& O\lp\frac{1}{t}\rp.
\end{align*}
Also,
\begin{align*}
  (\pt +i\px^2 ) \Phi_v^-=O\lp\frac{1}{t}\rp.
\end{align*}
We also define the asymptotic profiles
\begin{equation}\label{gamma pm defn}
    \begin{split}
         \gamma^+(t,v)=&\int \up (t,x) \widebar{\Phi_v^+}(t,x)\, dx\\
     \gamma^-(t,v)=&\int \um (t,x) \widebar{\Phi_v^-}(t,x)\, dx. 
    \end{split}
\end{equation}
More details about how these wave packets are chosen will be given in section \ref{Wave Packets section}.

\textbf{III. Bootstrap setup}
As in the case of DNLS, it will be necessary to use a bootstrap argument to complete the proof of the main theorem. We will make the bootstrap assumptions
\begin{align}
\label{bootstrap mbo}
    \|u^\pm\|_{L^\infty}\leq D\epsilon  \langle t\rangle^{-\frac{1}{2}}, \qquad  \|u_x^\pm\|_{L^\infty}\leq D\epsilon  \langle t\rangle^{-\frac{1}{2}},
\end{align}
where $D$ is a large universal constant that will be chosen later. Since $u=u^++u^-$ exactly, the same bounds hold for $u$ and $u_x$ with $D$ replaced by $2D$; these are the bounds used in the energy estimates of Section \ref{Vector field bnds section}.

\textbf{IV. Difference bounds}
Next, we need the equivalent of Lemma 3.1 of \cite{byars2024globaldynamicssmalldata}, comparing $\gamma^\pm$ with $u^\pm.$

\begin{lemma}\label{difference lemma mbo}
    For $u$ a solution to \eqref{mBO}, $\gamma^{\pm}$ as defined in \eqref{gamma pm defn}, and $c$ the constant of Lemma \ref{commutator lemma}, we have
    \begin{equation}
   \|\gamma^\pm\|_{L^\infty}\lesssim t^{\frac{1}{2}}\|u^\pm\|_{L^\infty}, \,\,\,\,\|\gamma^\pm\|_{L^2_v}\lesssim \|u^\pm\|_{L^2_x}, \,\,\,\,\|\gamma^\pm\|_{L^\infty_v}\lesssim t^{\frac14}\|u^\pm\|_{L^2_x}.
        \label{gamma+bnds}
    \end{equation}
Also,
\begin{align}\label{v gamma+ L^infty}
\|v\gamma^\pm\|_{L^\infty_v}\lesssim& t^{-\frac{3}{4}}\|Lu\|_{L^2_x}+t^{-\frac12}|c|+ t^{\frac{1}{2}}\|u_x^\pm\|_{L^\infty}.
\end{align}
We have the following spatial difference bounds,
\begin{align}
\|u^\pm(t,vt)-t^{-\frac{1}{2}}e^{i\phi^\pm(t,vt)}\gamma^\pm(t,v)\|_{L^\infty}\lesssim& t^{-\frac34}\|Lu\|_{L^2}+|c|t^{-\frac12}
      \label{spatial diff bnds mbo}
\end{align}
\begin{align}
    |v|\left|u^\pm(t,vt)-t^{-\frac{1}{2}}e^{i\phi^\pm(t,vt)}\gamma^\pm(t,v)\right|\lesssim & t^{-\frac34}\lp\|Lu\|_{H^1}+|c|\rp    \label{v spatial diff mbo}
\end{align}
spatial bounds for the derivative of our solution,
\begin{align}
\left\|u^\pm_x(t,vt)\mp \frac{i}{2}t^{-\frac{1}{2}}e^{i\phi^\pm(t,vt)}v\gamma^\pm(t,v)\right\|_{L^\infty}\lesssim& t^{-\frac34}\lp\|\px Lu\|_{L^2}+\|u\|_{L^2}\rp+ t^{-\frac54}\|Lu\|_{L^2}+t^{-1}|c|.
   \label{spatialu_xbnds+}
\end{align} 
\
\end{lemma}

The proof of this Lemma is very similar to the DNLS case and is given in Section \ref{diff lemma proof sec}.

\textbf{V. Asymptotic equation} Let $u$ be a solution to \eqref{mBO} and fix $\delta<\frac{1}{4}$. Then $\gp, \gm$ satisfy the following equations:
\begin{align} \label{gamma + eqn}
  \pt \gp =&  i\frac{v}{2t}(|\gp|^2\gp +2|\gm|^2 \gp)+R^+(t,v)\\
  \pt \gm =& -i\frac{v}{2t} \lp |\gm|^2\gm +2 |\gp|^2\gm \rp+R^-(t,v) \label{gamma - eqn}
\end{align}
where
\begin{equation}
    \begin{aligned}
        \label{R pm bounds}
    \|R^\pm\|_{L^\infty}\lesssim& t^{-\frac{5}{4}}\|Lu\|_{L^2} + |c|\,|v|\,t^{-\frac{1}{2}}\la t^{\frac{1}{2}}v\ra^{-N}+ D^2 \ep^3 t^{-\frac54+\delta}\\
    +& D^3 \ep^3 t^{-\frac{3}{2}+\delta}+D^2\ep^2 t^{-\frac{5}{4}}(\|Lu\|_{H^1}+|c|)
    \end{aligned}
\end{equation}
and
\begin{equation}
    \begin{aligned}
        \|vR^\pm\|_{L^\infty}\lesssim&D^3\ep^3 t^{-\frac32}+D^2\ep^3t^{-\frac{5}{4}}+ D^3 \ep^3 t^{-2+2\delta}+  |c|v^2t^{-\frac{1}{2}}\la t^{\frac{1}{2}}v\ra^{-N}\\
    &+  (D^2\ep^2  t^{-\frac54}+t^{-\frac{5}{4}})\lp\|Lu\|_{H^1}+|c|\rp+ D\ep t^{-\frac54} \lp\|Lu\|_{H^1}+|c|\rp^2.
    \end{aligned}
\end{equation}
In order to obtain these equations, a normal form-like analysis is needed, where we treat the resonant and nonresonant terms separately.
The proof will be given in section \ref{asymp eqn mbo sec}.

\textbf{VI. Vector field energy estimates}
As much of the proof relies on the vector field energy estimates, we will show how these can be obtained. Here, we prove:
\begin{lemma}
    \label{LuBound lemma mbo}
    Let $L$ be the operator defined in \eqref{L defn mbo}. Then, for any solution $u$ to \eqref{mBO} satisfying \eqref{initials mbo} and \eqref{bootstrap mbo}, we have
    \begin{align} \label{vector field EE}
    \|Lu\|_{L^2_x}\lesssim & \ \epsilon\langle t\, \rangle^{CD^2\epsilon^2},\\
        \|\px Lu\|_{L^2_x}\lesssim&  \ \epsilon\, \langle t\rangle^{CD^2\epsilon^2},\label{derivative vector field EE}
\end{align}
where $C$ is some universal constant which does not depend on $\ep$ or $D.$
\end{lemma}
The proof of this Lemma is given in Section \ref{Vector field bnds section}.

\textbf{VII. Closing the bootstrap argument}
Lastly, we need to improve the bootstrap bounds in order to close the bootstrap argument.  More precisely, we prove
\begin{lemma}
\label{close bootstrap lemma}
    Let $t\geq 1$ and $u$ be a solution to \eqref{mBO} satisfying \eqref{initials mbo} and \eqref{bootstrap mbo}. Then we have 
\begin{align*}
     \|u^\pm\|_{L^\infty}\leq E\epsilon  \langle t\rangle^{-\frac{1}{2}}, \qquad 
          \|u_x^\pm\|_{L^\infty}\leq F\epsilon  \langle t\rangle^{-\frac{1}{2}},
\end{align*}
    where $E,F<D$.
\end{lemma}
This closes the bootstrap argument and completes the proof of Theorem \ref{mbo thm}.

\subsection{Positive and Negative Frequency Projection norms}

Since we are projecting to positive and negative frequencies, it's important to see how the $L^p$ norms of $u^\pm$ are related to those of $u$.
In general, for some function $f,$ let $\fp, \fm$ be the projections to positive and negative frequencies, as given in Definition \ref{LP defn}, so that $f=f^++f^-$. Then,
\begin{align*}
    \|f\|_{L^2_x}^2
=&\|\fp\|_{L^2_x}^2+\|\fm\|_{L^2_x}^2.
\end{align*}
In particular, if $\|u_0\|_{L^2}\leq \ep,$ then also $\|\up_0\|_{L^2}, \|\um_0\|_{L^2}\leq \ep.$

\section{Wave packets \label{Wave Packets section}}
We now proceed using the method of testing by wave packets, which was developed by Ifrim and Tataru for NLS \cite{cite:NLS}, and explained in general in \cite{cite:WavePax}. This method has been used numerous times to prove low regularity well-posedness, dispersive decay results, and modified scattering \cite{albert.ovidiu2023gsqg,albert.ovidiu2024sqg,aiifrimtataru2022two,byars2024globaldynamicssmalldata, byars2026modscatteringAC, cite:NLS}.  The idea is to test the solution against a wave packet which is localized in space along rays $x=vt$, and in frequency (with respect to the uncertainty principle).  

In this paper, we will adapt the method to deal with the positive and negative frequency projections of the equation separately. First, we give the formal derivation of the wave packets used in the proof.  In \cite{albert.ovidiu2023gsqg}, they give a general form of the wave packets, based on \cite{cite:WavePax}. For a general dispersive equation
\begin{align*}
    (i\pt -A)u=0
\end{align*}
where $a$ is the symbol of $A,$ i.e. the dispersion relation, define $v=a'(\xi),$ and $\phi=t(v\xi_v-a(\xi_v))$.  
Then the wave packets take the form
\begin{align*}
    \Phi_v=|a''(\xi_v)|^{-\frac{1}{2}}e^{i\phi}\chi\lp \frac{x-tv}{t^{\frac{1}{2}}|a''(\xi_v)|^{\frac{1}{2}}}\rp.
\end{align*}
For \eqref{mBO}, we have the distributional derivatives
\begin{align*}
    a(\xi)=& \xi|\xi|\\
    a'(\xi)=&2|\xi|\\
    a''(\xi)=&\begin{cases}
        2, & \xi>0\\
        -2, & \xi<0.
    \end{cases}
\end{align*}
This is yet another suggestion that splitting by positive and negative frequencies is a good idea.

 \subsection{Positive frequency wave packets\label{positive wave paxs}}

Consider the linear equation projected to positive frequencies
\begin{align*}
    \pt u^+-i\px^2 u^+=0
\end{align*}
Here we have
$v=2\xi$, so we set $\xi_v=\frac{1}{2}v.$
For the positive frequency case, we have the phase
\[\phi^+=t\lp \frac{v^2}{2}-\xi_v^2\rp=t\lp\frac{v^2}{2}-\lp \frac{v}{2}\rp^2\rp=\frac{tv^2}{4}.\]
For $v=x/t,$ this is $\phi^+=\frac{x^2}{4t}.$ Therefore the wave packets are given by
\begin{align*}
    \Phi_v^+(t,x)=e^{i\frac{x^2}{4t}}\chi\lp \frac{x-tv}{t^{\frac{1}{2}}}\rp,
\end{align*}
where $\chi$ is a real, even bump function with $\int \chi(y)\, dy =1.$ The wave packets are approximate solutions to the positive frequency linear equation:
\begin{align*}
    \pt \Phi_v^+-i\px^2\Phi_v^+    =&\frac{e^{i\frac{x^2}{4t}}}{2t} \px \lb (x-tv)\chi\lp \frac{x-tv}{t^{\frac{1}{2}}}\rp -2it^{\frac{1}{2}} \chi'\lp \frac{x-tv}{t^{\frac{1}{2}}}\rp  \rb\\
    =& O_{L^\infty}\lp\frac{1}{t}\rp.
\end{align*}

\subsection{Negative frequency wave packets}
We can do the same thing for the negative frequencies. The negative frequency projection linear equation is
\begin{align*}
    \pt u^-+i\px^2 u^-=0.
\end{align*}
In this case,
 $v=-2\xi$, so $ \xi_v =-\frac{1}{2}v.$
Then
$\phi^- =-\frac{x^2}{4t}$, and the wave packets associated to $\um$ are
\begin{align*}
    \Phi_v^-(t,x)=e^{-i\frac{x^2}{4t}}\chi\lp \frac{x-tv}{t^{\frac{1}{2}}}\rp.
\end{align*}
Similarly, we have
\begin{align*}
      (\pt +i\px^2 ) \Phi_v^-=O\lp\frac{1}{t}\rp.
\end{align*}
We then define asymptotic profiles $\g^{\pm}$ as
\begin{align*}
    \gamma^+(t,v)=&\int \up (t,x) \widebar{\Phi_v^+}(t,x)\, dx\\
     \gamma^-(t,v)=&\int \um (t,x) \widebar{\Phi_v^-}(t,x)\, dx.
\end{align*}

\subsection{\texorpdfstring{Preliminary facts about $\gp, \gm$}{Preliminary facts about gamma} }
We now make some observations about $\gp, \gm$.  First, define
\begin{equation*}
    w^+:=e^{-i\phi^+}\up (t,x), \,\,\,\, \wm:= e^{-i\phi^-}\um (t,x)
\end{equation*}
Then, we can write $\gp$ as a convolution

\begin{align*}
  \gamma^+=&\int \up(t,x) \widebar{\Phi_v^+}(t,x)\, dx=\int   \up(t,x) e^{-i\frac{x^2}{4t}}\chi\lp \frac{x-tv}{t^{\frac{1}{2}}}\rp\, dx\\
  =& \int \wp(t,x)\chi\lp \frac{x-tv}{t^{\frac{1}{2}}}\rp\, dx\\
  =&\int t \wp(t,ty)\chi\lp t^{\frac{1}{2}}(y-v)\rp\, dy\\
  =&t^{\frac{1}{2}}\wp(t,vt) * t^{\frac{1}{2}}\chi(t^{\frac{1}{2}}v).
\end{align*}
Similarly,
\begin{align*}
    \gm = t^{\frac{1}{2}}\wm(t,vt) * t^{\frac{1}{2}}\chi(t^{\frac{1}{2}}v).
\end{align*}
We also have
\begin{align*}
\partial_v[w^\pm(t,tv)]=t\partial_xw^\pm(t,tv)=\mp\frac{i}{2}e^{-i\phi^\pm}Lu^\pm(t,tv).
\end{align*}
where
\begin{align*}
    Lu^\pm = xu^\pm\pm 2ti\px u^\pm.
\end{align*}

\section{Proof of Lemma \ref{difference lemma mbo} \label{diff lemma proof sec}}
Let $j^\pm:=e^{-i\phi^\pm}u^\pm_x$ and let $\supp\chi\subset[-C_\chi,C_\chi]$. We first prove a Lemma that will help prove Lemma  \ref{difference lemma mbo}.
\begin{lemma}[Local differences]\label{local diff lemma}
Let $t\geq1$, $v\in\R$, and $|x-vt|\leq C_\chi t^{\frac12}$. Then
\begin{align}
    |w^\pm(t,x)-w^\pm(t,vt)|\lesssim& t^{-\frac34}\|Lu\|_{L^2}+|c|t^{-\frac12},\label{local diff 1}\\
   \,|v w^\pm(t,x)-v w^\pm(t,vt)|\lesssim& t^{-\frac34}\lp\|Lu\|_{H^1}+|c|\rp,\label{local diff 2}\\
    |j^\pm(t,x)-j^\pm(t,vt)|\lesssim&\ t^{-\frac34}\lp\|\px Lu\|_{L^2}+\|u\|_{L^2}\rp.\label{local diff 3}
\end{align}
\end{lemma}
\begin{proof}
We show the proof for the $+$ case and note that the $-$ case is identical. Recall from the above identities and Lemma \ref{commutator lemma}  we have
\begin{align}\label{dw}
    \px w^+=-\frac{i}{2t}e^{-i\phi^+}\lp P^+Lu+c\rp,\qquad \px j^+=-\frac{i}{2t}e^{-i\phi^+}P^+(\px Lu-u),
\end{align}
where for the second identity we used $L(u^+_x)=P^+L(u_x)$ and $Lu_x=\px Lu-u$. Then, writing $I_x$ for the interval between $vt$ and $x$, which has length at most $C_\chi t^{\frac12}$,
\begin{align*}
    |w^+(x)-w^+(vt)|\leq&\int_{I_x}|\px w^+|\leq\frac{1}{2t}\|P^+Lu\|_{L^2}|I_x|^{\frac12}+\frac{|c|}{2t}\Big|\int_{I_x}e^{-i\frac{y^2}{4t}}dy\Big| \\
    \lesssim&t^{-\frac34}\|Lu\|_{L^2}+\frac{|c|}{2t}\Big|\int_{I_x}e^{-i\frac{y^2}{4t}}dy\Big| 
\end{align*}
For the second we use two bounds on the Fresnel integral.  In general, we have for $a>0,$
\begin{align*}
    \int_{-\infty}^\infty e^{\pm ia x^2}\, dx=\sqrt{\frac{\pi}{a}}e^{\pm i\pi /4},
\end{align*}
which gives
\begin{align*}
    \left|\int_a^be^{-i\frac{y^2}{4t}}dy\right|\leq C_0t^{\frac12} \text{ for all } a<b.
\end{align*}
This proves \eqref{local diff 1}.
For \eqref{local diff 2}, 
\begin{align*}
    |vw^+(x)-vw^+(vt)|\leq\left|\int_{I_x}v\px w^+\right|\leq  \Big|\underbrace{\int_{I_x}  \frac{y}{t}\px w^+(t,y)\, dy}_{A} \Big|+\Big|\underbrace{\int_{I_x} \frac{y-vt}{t}\px w^+(t,y) \, dy}_B\Big|.
\end{align*}
Then,
\begin{align*}
    A=&\int_{I_x}   \frac{y}{t}\px (\wp (t,y))\, dy =\int_{I_x}   \frac{y}{t}(-\frac{i}{2t}e^{-i\phi^+}\lp P^+Lu+c\rp)\, dy\\
    =&\int_{I_x}  \frac{1}{t}(\p_y e^{-i\frac{y^2}{4t}}\lp P^+Lu+c\rp)\, dy\\
    =&\frac{1}{t}e^{-i\frac{y^2}{4t}}\lp P^+Lu+c\rp\Big|_{I_x}-\frac{1}{t}\int_{I_x}e^{-i\frac{y^2}{4t}}\p_y\lp P^+Lu+c\rp) \, dy
\end{align*}
By Sobolev embedding, this gives
\begin{align*}
    |A|\leq& \frac{1}{t}\|P^+ Lu\|_{L^\infty}+ \frac{|c|}{t}+\frac{1}{t}|I_x|^{\frac{1}{2}}\|\px P^+ Lu\|_{L^2}\\
    \lesssim& \frac{1}{t}\| Lu\|_{H^1}+ \frac{|c|}{t}+t^{-\frac{3}{4}}\|\px Lu\|_{L^2}
\end{align*}
This proves the bound for $A.$  For $B,$ first note that if $y\in I_x,$ then $\frac{|y-vt|}{t}\leq C_\chi t^{-\frac{1}{2}}.$ Then,
\begin{align*}
|B|=& \Big| \int_{I_x} \frac{y-vt}{t}\px w^+(t,y) \, dy \Big| \lesssim t^{-\frac{1}{2}}\int_{I_x}|\px w^+|\\
\leq& \frac{1}{2t^{\frac{3}{2}}}\|P^+Lu\|_{L^2}|I_x|^{\frac12}+\frac{|c|}{2t^{\frac{3}{2}}}\Big|\int_{I_x}e^{-i\frac{y^2}{4t}}dy\Big|\\
\lesssim& t^{-\frac{5}{4}}\|Lu\|_{L^2}+|c|t^{-1}.
\end{align*}
 Finally, \eqref{local diff 3} follows from the second identity in \eqref{dw} exactly as the first term in \eqref{local diff 1}.
\end{proof}

\begin{proof}[Proof of Lemma \ref{difference lemma mbo}]
Define $\chi_{v,t}(x):=\chi\lp\frac{x-vt}{t^{\frac{1}{2}}}\rp$, and recall the identities
\begin{align*}
w^{\pm}=&e^{-i\phi^{\pm}} u^{\pm}, \hspace{2.05cm}   \pv(w^{\pm})= \mp \frac{i}{2}e^{-i\phi^{\pm}}Lu^{\pm}\\
    \gamma^{\pm}=&\int u^{\pm}\overline{\Phi_v^{\pm}(t,x)}\, dx,  \hspace{1.3cm}    \gamma^{\pm} =t^{\frac{1}{2}}w^{\pm}(t,vt) * t^{\frac{1}{2}}\chi(t^{\frac{1}{2}}v).
\end{align*}
Note that 
\begin{align} \label{chi norms}
    \|\chi_{v,t}\|_{L^1}\lesssim t^{\frac{1}{2}}, \,\,\,   \|\chi_{v,t}\|_{L^2}\lesssim t^{\frac{1}{4}}, \,\,\,    \|\px \chi_{v,t}\|_{L^1}\lesssim 1.
\end{align}
Since $\chi$ is real, 
\[\gamma^\pm(t,v)=\int u^\pm\overline{\Phi_v^\pm}\,dx=\int w^\pm(t,x)\chi_{v,t}(x)\,dx.\]
Then the first inequality of \eqref{gamma+bnds} is
\begin{align*}
\|\gamma^\pm\|_{L^\infty_v}\lesssim& \|\chi_{v,t}\|_{L^1}\|w^\pm\|_{L^\infty}
\leq t^{\frac{1}{2}} \|\chi\|_{L^1}\|u^\pm\|_{L^\infty}\lesssim t^{\frac{1}{2}} \|u^\pm\|_{L^\infty}.
\end{align*}
We obtain the second inequality by applying Young's inequality to $\g^\pm$ written as a convolution to get
\begin{align*}
    \|\gamma^\pm\|_{L^2_v}\leq t\|w^\pm(t,t\cdot)\|_{L^2}\|\chi(t^{\frac12}\cdot)\|_{L^1}=\|u^\pm\|_{L^2_x}\|\chi\|_{L^1_x}\lesssim \|u^\pm\|_{L^2_x}.
\end{align*}
For the third one, we have
\begin{align*}
\|\gamma^\pm\|_{L^\infty_v}\lesssim& \|\chi_{v,t}\|_{L^2}\|w^\pm\|_{L^2}
\leq t^{\frac{1}{4}} \|u^\pm\|_{L^2}.
\end{align*}
Next, for the difference bounds,
\begin{align*}
   |u^\pm(t,vt)-t^{-\frac12}e^{i\phi^\pm(t,vt)}\gamma^\pm|=& |w^\pm(t,vt)-t^{-\frac12}\gamma^\pm(t,v)|\\
   =&\left|t^{-\frac12}\int\lp w^\pm(t,x)-w^\pm(t,vt)\rp\chi_{v,t}(x)\,dx \right|
\end{align*}
Then, the bounds \eqref{spatial diff bnds mbo} and \eqref{v spatial diff mbo} follow by applying $t^{-\frac12}\int\chi_{v,t}\,dx=1$ along with \eqref{local diff 1} and \eqref{local diff 2} respectively.

For \eqref{v gamma+ L^infty} and \eqref{spatialu_xbnds+}, we use the localization $v\approx x/t$ along with the definition of the vector field $L$, and Lemma \ref{commutator lemma}
\begin{align*}
    v\gamma^{\pm}(t,v)=&\int v u^{\pm}\overline{\Phi_v^{\pm}(t,x)}\, dx=\int \lp \frac{x}{t}-\frac{x-vt}{t} \rp u^{\pm}\overline{\Phi_v^{\pm}(t,x)}\, dx\\
    =& \int \frac{1}{t}\lp Lu^\pm \mp 2ti\px  u^\pm\rp\overline{\Phi_v^{\pm}(t,x)}\, dx- \int \frac{x-vt}{t} u^{\pm}\overline{\Phi_v^{\pm}(t,x)}\, dx\\
    =& \underbrace{\int \left[\frac{P^\pm Lu\pm c}{t}\mp 2iu^{\pm}_x\right]e^{-i\phi^{\pm}}\chi_{v,t}\, dx}_{I}\underbrace{- \int \frac{x-vt}{t} w^\pm\,\chi_{v,t}\,dx}_{II}
\end{align*}
The first term is bounded by 
\begin{align*}
|I|\lesssim \frac1t\lp\|Lu\|_{L^2}\|\chi_{v,t}\|_{L^2}+|c|\|\chi_{v,t}\|_{L^1}\rp+2\|u^\pm_x\|_{L^\infty}\|\chi_{v,t}\|_{L^1}\lesssim t^{-\frac34}\|Lu\|_{L^2}+t^{-\frac12}|c|+t^{\frac12}\|u^\pm_x\|_{L^\infty}.
\end{align*}
By a change of variables, the second is bounded by
\begin{align*}
|II|\lesssim    \left|\int \frac{x-vt}{t} w^\pm\,\chi_{v,t}\,dx\right|\lesssim \|w^\pm\|_{L^\infty}\int |y|\chi(y)\, dy\lesssim \|u^\pm\|_{L^\infty}.
\end{align*}
However, this does not have the right form. Instead, we will rewrite the term as 
\begin{align*}
    II=-t^{-\frac{1}{2}}\int w^\pm\,(y\chi)_{v,t}\,dx .
\end{align*}
Now define $\Theta(y):=\int_{-\infty}^y s\chi(s)\,ds$.  Since $\chi$ is even, $\int y\chi(y)\,dy=0$.  Therefore, $\Theta(y)$ is again a bump function, with $\Theta'(y)=y\chi(y)$.  Since $\px[\Theta_{v,t}]=t^{-\frac12}(y\chi)_{v,t}$, we may integrate by parts and use \eqref{dw} to get
\begin{align*}
  |II|=&\left|\int \px w^\pm\,\Theta_{v,t}\,dx\right|=\frac{1}{2t}\left|\int e^{-i\phi^\pm}\lp P^\pm Lu\pm c\rp\Theta_{v,t}\,dx\right|,\\
    |II|\leq& \frac{1}{2t}\lp\|Lu\|_{L^2}\|\Theta_{v,t}\|_{L^2}+|c|\|\Theta_{v,t}\|_{L^1}\rp\lesssim t^{-\frac34}\|Lu\|_{L^2}+t^{-\frac{1}{2}}|c|.
\end{align*}
This proves \eqref{v gamma+ L^infty}. Lastly, we prove \eqref{spatialu_xbnds+}.  From before we had
\begin{align*}
    v\gamma^{\pm}(t,v)
    =& \int \left[\frac{P^\pm Lu\pm c}{t}\mp 2iu^{\pm}_x\right]e^{-i\phi^{\pm}}\chi_{v,t}\, dx- \int \frac{x-vt}{t} w^\pm\,\chi_{v,t}\,dx.
\end{align*}
Multiplying both sides by $\mp \frac{i}{2} t^{-\frac12}$ gives
\begin{align*}
   \mp \frac{i}{2} t^{-\frac12} v\gamma^{\pm}(t,v)
    =& \mp \frac{i}{2} t^{-\frac32} \int \left[P^\pm Lu\pm c\right]e^{-i\phi^{\pm}}\chi_{v,t}\, dx - t^{-\frac12} \int  j^\pm(x) \chi_{v,t}\, dx \\
    &\pm \frac{i}{2} t^{-\frac12} \int \frac{x-vt}{t} w^\pm\,\chi_{v,t}\,dx.
\end{align*}
Then adding $j^\pm(t,vt)$ to both sides results in
\begin{align*}
   j^\pm(t,vt)\mp \frac{i}{2} t^{-\frac12} v\gamma^{\pm}(t,v)
    =&  - t^{-\frac12}\int(j^\pm(x)-j^\pm(vt))\chi_{v,t}\,dx\\
    &\mp \frac{i}{2} t^{-\frac32} \int \left[P^\pm Lu\pm c\right]e^{-i\phi^{\pm}}\chi_{v,t}\, dx \underbrace{\pm \frac{i}{2} \int \frac{x-vt}{t^{\frac{3}{2}}} w^\pm\,\chi_{v,t}\,dx}_{III}.
\end{align*}
Taking the absolute value,
\begin{align*}
    \Big| j^\pm(t,vt)\mp \frac{i}{2} t^{-\frac12} v\gamma^{\pm}(t,v)\Big|\lesssim \sup_{|x-vt|\leq C_\chi t^{\frac{1}{2}}}|j^\pm(x)-j^\pm(vt)|\,\|\chi\|_{L^1}+ t^{-\frac54}\|Lu\|_{L^2}+t^{-1}|c|+ |III|.
\end{align*}
By applying \eqref{local diff 3}, we have
\begin{align*}
\left\|u^\pm_x(t,vt)\mp \frac{i}{2}t^{-\frac{1}{2}}e^{i\phi^\pm(t,vt)}v\gamma^\pm(t,v)\right\|_{L^\infty}\lesssim& t^{-\frac34}\lp\|\px Lu\|_{L^2}+\|u\|_{L^2}\rp+ t^{-\frac54}\|Lu\|_{L^2}+t^{-1}|c|+ |III|.
\end{align*} 
To bound $|III|,$ we do the same procedure as for $II$ above.
\begin{align*}
    III=\pm \frac{i}{2}t^{-1}\int w^\pm\,(y\chi)_{v,t}\,dx .
\end{align*}
For $\Theta(y):=\int_{-\infty}^y s\chi(s)\,ds$ as defined before, we integrate by parts to get
\begin{align*}
  III=&\mp \frac{i}{2}t^{-\frac{1}{2}}\int \px w^\pm\,\Theta_{v,t}\,dx=\mp\frac{1}{4t^{\frac{3}{2}}}\int e^{-i\phi^\pm}\lp P^\pm Lu\pm c\rp\Theta_{v,t}\,dx,\\
    |III|\leq& \frac{1}{2t^{\frac{3}{2}}}\lp\|Lu\|_{L^2}\|\Theta_{v,t}\|_{L^2}+|c|\|\Theta_{v,t}\|_{L^1}\rp\lesssim t^{-\frac54}\|Lu\|_{L^2}+t^{-1}|c|.
\end{align*}

\end{proof}

\section{The Bootstrap Argument \label{bootstrap section}}
In order to complete the proof, it will be necessary to use a bootstrap argument. As a reminder, the assumptions we make here are
\begin{align}\label{bootstrap equation}
    \|u^\pm\|_{L^\infty}\leq D\epsilon  \langle t\rangle^{-\frac{1}{2}}, \qquad  \|u_x^\pm\|_{L^\infty}\leq D\epsilon  \langle t\rangle^{-\frac{1}{2}},
\end{align}
where $D$ is a large universal constant that will be chosen later. The goal will then be to prove all of these bounds but with a constant strictly less than $D.$  This will be possible assuming $\epsilon$ is sufficiently small.

\begin{remark}
 In order to do this bootstrap argument, we require the map $t \mapsto \|u_x^\pm(t)\|_{L^\infty}$ to be continuous.  However, this does not hold for general $H^1$ solutions. Therefore, we will prove the estimates assuming smooth initial data $u_0 \in H^s(\mathbb{R})$ ($s \gg 1$). The result for general $H^1$ initial data then follows by a standard density and approximation argument.

\end{remark}

\section{\texorpdfstring{$\gamma^\pm$ lemma}{gamma lemma} \label{gamma lemma section}}

\subsection{\texorpdfstring{Asymptotic equation for $\gamma^\pm$}{Asymptotic equation for gamma}\label{asymp eqn mbo sec} }

In this section, we will derive the asymptotic equation for $\gamma^{\pm}.$  We do this by writing out the derivative of $\gamma^+$, and then term-by-term extracting the principal part of the equation.
We compute
\begin{align*}
    \pt \gamma^+=&\int  \pt \up (t,x) \widebar{\Phi_v^+}(t,x)+\up (t,x) \pt \widebar{\Phi_v^+}(t,x)\, dx\\
    =&\int [i\px^2 \up +P^+\px (|u|^2u)] \widebar{\Phi_v^+}(t,x)+\up (t,x) \pt \widebar{\Phi_v^+}(t,x)\, dx\\
    =&\int [\pt \widebar{\Phi_v^+}+i\px^2 \widebar{\Phi_v^+}] \up + P^+\px (|u|^2 u) \widebar{\Phi_v^+}\, dx\\
    =& \underbrace{\int \overline{\frac{e^{i\frac{x^2}{4t}}}{2t} \px \lb (x-tv)\chi\lp \frac{x-tv}{t^{\frac{1}{2}}}\rp -2it^{\frac{1}{2}} \chi'\lp \frac{x-tv}{t^{\frac{1}{2}}}\rp  \rb} \up \, dx}_{R_1^+} - \int P^+(|u|^2 u) \px\widebar{\Phi_v^+}\, dx.
\end{align*}
The first term is bounded as a part of the remainder. For the second, we first note that $\px \Phi_v^+=\frac{ix}{2t}\Phi_v^++t^{-\frac{1}{2}}\Psi_v^+$, where $\Psi_v^+:=e^{i\phi^+}\chi'_{v,t},$ which gives
\begin{align*}
    - \int P^+(|u|^2 u) \px\widebar{\Phi_v^+}\, dx= \underbrace{\int  \frac{ix}{2t} P^+(|u|^2 u)\widebar{\Phi_v^+}\, dx}_{A}\underbrace{-\int t^{-\frac{1}{2}}\widebar{\Psi_v^+} P^+(|u|^2 u)\, dx}_{R_2^+}.
\end{align*}
The second term, $R_2^+,$ has enough decay so it is the second remainder term. For the first, we need to decompose the nonlinearity into its positive and negative frequency parts to extract the principal terms.  By writing $u=\up+\um$ and expanding, we have
\begin{align*}
    |u|^2u=&|\up|^2\up+(\up)^2\widebar{\um}+2|\up|^2\um + 2\up |\um|^2+(\um)^2\widebar{\up}+|\um|^2\um\\
P^+(|u|^2u)=&P^+(|\up|^2\up)+(\up)^2\widebar{\um}+2P^+(|\up|^2\um )+ 2P^+(\up |\um|^2)+P^+(|\um|^2\um),
\end{align*}
where $P^+((\up)^2\widebar{\um})=(\up)^2\widebar{\um}$  because each factor is supported on positive frequencies, and  $P^+((\um)^2\widebar{\up})=0$ since each factor is supported on negative frequencies. 
We then write 
\begin{align*}
 A=&A_1+A_2+A_3+A_4+A_5,\qquad A_j:=\frac{i}{2t}\int x\,P^+(f_j)\,\widebar{\Phi_v^+}\,dx,\\
 f_1:=&|\up|^2\up,\quad f_2:=(\up)^2\widebar{\um},\quad f_3:=2|\up|^2\um,\quad f_4:=2\up|\um|^2,\quad f_5:=|\um|^2\um .
\end{align*}

On a first attempt, one might compare each term in the expansion to its corresponding $\gamma^\pm$. This would give
\begin{align}\label{gamma + eqn preliminary}
    \pt \gp =&  i\frac{v}{2t}(|\gp|^2\gp + e^{\frac{iv^2t}{2}}(\gp)^2 \widebar{\gm}+2e^{\frac{-iv^2t}{2}}\gm |\gp|^2 +2|\gm|^2 \gp)+r^+
\end{align}
where
\begin{align*}
    r^+=& R_1^++R_2^++r_3^++r_4^+,
\end{align*}
and
\begin{align*}
    r_3^+=&\int \frac{ix}{2t} P^+(|\up|^2\up)\widebar{\Phi_v^+}\, dx - \frac{iv}{2}\int \up \widebar{P^+\Phi_v^+} |\up(t,vt)|^2 \,dx \\
    &+\int \frac{ix}{2t}  P^+((\up)^2\widebar{\um})\widebar{\Phi_v^+}\, dx- \frac{iv}{2} \int \up(t,vt)\widebar{\um(t,vt)} \up(t,x) \widebar{P^+\Phi_v^+}(t,x)\, dx\\
    &+ \int  \frac{ix}{t} P^+(|\up(t,x)|^2 \um(t,x))\widebar{\Phi_v^+} \, dx -iv\int \um(t,vt)\widebar{\up(t,vt)} \up(t,x) \widebar{P^+\Phi_v^+}(t,x)\, dx\\
    &+ \int \frac{ix}{t}  P^+(\up |\um(t,x)|^2) \widebar{\Phi_v^+}  \, dx-i v\int  |\um(t,vt)|^2\up(t,x) \widebar{P^+\Phi_v^+} \,dx \\
&+\int \frac{ix}{2t} P^+(|\um|^2\um)\widebar{\Phi_v^+}\, dx -\frac{iv}{2}\int  |\um(t,vt)|^2 \um(t,x) \widebar{P^+\Phi_v^+} \,dx 
\end{align*}
\begin{align*}
   r_4^+=& i\frac{v}{2}\Big[\lp |\up(t,vt)|^2-t^{-1}|\gp|^2\rp+\lp \up(t,vt) \widebar{\um}(t,vt)- t^{-1}e^{\frac{iv^2t}{2}}\gp \widebar{\gm}\rp\\
   +& 2\lp\um(t,vt)\widebar{\up}(t,vt)-t^{-1}e^{\frac{-iv^2t}{2}}\gm \widebar{\gp}\rp+2\lp |\um(t,vt)|^2-t^{-1}|\gm|^2 \rp \Big]\gp(t,v).
\end{align*}

In equation \eqref{gamma + eqn preliminary}, the first and fourth terms are called resonant terms. When the energy estimate is performed, these terms become zero. 
However, the second and third terms in this equation are nonresonant, and when performing the energy estimate, they do not have enough decay. 

At this point, there are two options. The first option is to do a normal form transformation to transform the middle two terms away.  The other option is to notice that $A_j$ for $j=2,3,5$ can be bounded directly.  This will be simpler, so we will opt for this route. The way these terms are bounded is by using the exponential to write $\frac{x}{2t}e^{i\frac{x^2}{4t}}$ as a derivative, then applying integration by parts.  The idea is given in Lemma \ref{nonres lemma}. The reason this cannot be done for $A_1, A_4$ is because they do not have an exponential factor.

To get this version of the equation, we return to where we had
\begin{align*}
    \pt \gamma^+= \sum_{j=1}^5 A_j+R_1^++R_2^+
\end{align*}
$A_1$ and $A_4$ are the resonant terms, so we extract the principal terms of the equation from those. For $A_1,$
\begin{align*}
   A_1=& \int \frac{ix}{2t} P^+(|\up|^2\up)\widebar{\Phi_v^+}\, dx\\
     =&\underbrace{ \int \frac{ix}{2t} P^+(|\up|^2\up)\widebar{\Phi_v^+}\, dx - \frac{iv}{2}\int \up \widebar{P^+\Phi_v^+} |\up(t,vt)|^2 \,dx }_{R_{3,1}^+} + \underbrace{i\frac{v}{2}|\up(t,vt)|^2\gamma^+(t,v)}_{B_1}.
\end{align*} 

For $A_4,$
\begin{align*}
     A_4 =&\underbrace{\int \frac{ix}{t}  P^+(\up |\um(t,x)|^2) \widebar{\Phi_v^+}  \, dx- \int i v |\um(t,vt)|^2\up(t,x) \widebar{P^+\Phi_v^+} \,dx }_{R_{3,4}^+} + \underbrace{iv|\um(t,vt)|^2\gamma^+(t,v)}_{B_4}.
\end{align*}
The nonresonant terms become part of the remainder, so define $R_{3,j}^+:=A_j$ for $j=2,3,5.$ 
Next, each $\up(t,vt)$ is compared to $t^{-\frac{1}{2}}e^{i\phi^+(t,vt)}\gp(t,v) =t^{-\frac{1}{2}}e^{\frac{iv^2t}{4}}\gp(t,v) $ and each $\um(t,vt)$ is compared to $t^{-\frac{1}{2}}e^{i\phi^-(t,vt)}\gm(t,v)  =t^{-\frac{1}{2}}e^{\frac{-iv^2t}{4}}\gm(t,v) $.
This means that we must add and subtract the term:
\begin{align*}
 i\frac{v}{2t}(|\gp|^2+2|\gm|^2 )\gp(t,v).
\end{align*}
This gives the equation:
\begin{align}
    \pt \gp =&  i\frac{v}{2t}(|\gp|^2\gp  +2|\gm|^2 \gp)+R^+,
\end{align}
where 
\begin{align*}
    R^+=& R_1^++R_2^++R_3^++R_4^+,
\end{align*}
and
\begin{align*}
   R_4^+=& i\frac{v}{2}\Big[\lp |\up(t,vt)|^2-t^{-1}|\gp|^2\rp+2\lp |\um(t,vt)|^2-t^{-1}|\gm|^2 \rp \Big]\gp(t,v).
\end{align*}

\subsection{Bounding the remainders \label{R bnds mbo sec}}

In this section, we will bound the remainder terms.  In order to do this, we first need to look at a few common quantities. 

\begin{lemma}\label{Packet tails lemma} Let $t\geq 1,  \Phi_v^\pm (t,x)=e^{i\phi^\pm}\chi\lp \frac{x-vt}{\sqrt{t}}\rp,\,\Psi_v^\pm:=e^{i\phi^\pm}\chi'\lp \frac{x-vt}{\sqrt{t}}\rp$. Then,
    \begin{align*}
    \|\Phi_v^+\|_{L^2}\approx& t^{\frac{1}{4}}, \quad      \|P^+\Phi_v^+\|_{L^2}\lesssim t^{\frac{1}{4}}\\
    \|\Psi_v^+\|_{L^2}\approx& t^{\frac{1}{4}}, \quad  \|P^+\Psi_v^+\|_{L^2}\lesssim t^{\frac{1}{4}}.
\end{align*}
We also have packet tails in certain cases.  For $v> 0,$
\begin{align}\label{P+ L^2 bound}
    \|P^+ \overline{\Phi_v^+}\|_{L^2}&=  \|P^- \Phi_v^+\|_{L^2}  \lesssim t^{\frac{1}{4}}\la t^{\frac{1}{2}}v\ra^{-N}, \,\, \forall N\geq 0,\\
   \|P^+(x\widebar{\Phi_v^+})\|_{L^2}&\lesssim_N t^{\frac34}\la t^{\frac12}v\ra^{1-N},\,\, \forall N\geq 0.
\end{align}
When $v<0,$ we have
\begin{align}\label{P- L^2 bound}
    \|P^- \overline{\Phi_v^+}\|_{L^2} &=    \|P^+ \Phi_v^+\|_{L^2}\lesssim t^{\frac{1}{4}}\la t^{\frac{1}{2}}v\ra^{-N}, \,\, \forall N\geq 0\\
    \| P^- (x \overline{\Phi_v^+}) \|_{L^2} &\lesssim_N t^{\frac{3}{4}} \la t^{\frac{1}{2}} v \ra^{1-N}, \,\, \forall N\geq 0.
\end{align}
\end{lemma}

\begin{proof}
First, 
\begin{align*}
     \|\Phi_v^+\|_{L^2}= \|\chi_{v,t}\|_{L^2}\approx t^{\frac{1}{4}}, 
\end{align*}
and 
\begin{align*}
    \|P^+\Phi_v^+\|_{L^2}\leq \|\Phi_v^+\|_{L^2}.
\end{align*}
The bounds for $\Psi_v^+$ follow similarly. 
For the tails, consider the Fourier transform of $P^+\Phi_v^+.$ First,
\[
\big(\widebar{\Phi_v^+}\big)^{\wedge}(\xi) = t^{\frac{1}{2}} e^{-i t \left(\frac{v^2}{4} + v\xi\right)} \widehat{g}\left( t^{\frac{1}{2}} \left(\xi + \frac{v}{2}\right) \right)
\]
where $g(y) = e^{-i \frac{y^2}{4}} \overline{\chi}(y)$, so $\widehat{g}$ is Schwartz.
Then, 
\[
\|P^+ \overline{\Phi_v^+}\|_{L^2} = {t^{\frac{1}{4}}}\left( \int_{\frac{1}{2}t^{\frac{1}{2}}v}^{\infty} |\widehat{g}(\eta)|^2 d\eta \right)^{\frac{1}{2}},
\]
and using that $\wh{g}$ is Schwartz, we get that for $v>0,$
\begin{align*}
    \|P^+ \overline{\Phi_v^+}\|_{L^2} \lesssim t^{\frac{1}{4}}\la t^{\frac{1}{2}}v\ra^{-N}, \,\, \forall N\geq 0.
\end{align*}
On the other hand, when $v<0$ we have
\begin{align*}
    \|P^- \overline{\Phi_v^+}\|_{L^2} \lesssim t^{\frac{1}{4}}\la t^{\frac{1}{2}}v\ra^{-N}, \,\, \forall N\geq 0.
\end{align*}
For the $\|P^+(x\widebar{\Phi_v^+})\|_{L^2}$ bound, using the property $\widehat{x f}(\xi) = i \partial_\xi \widehat{f}(\xi)$, we differentiate with respect to $\xi$:
\begin{align*}
\widehat{x \overline{\Phi_v^+}}(\xi) &= i \partial_\xi \left[ t^{\frac{1}{2}} e^{-i t \left( \frac{v^2}{4} + v \xi \right)} \widehat{g}\left( t^{\frac{1}{2}} \left( \xi + \frac{v}{2} \right) \right) \right] \\
&= t^{\frac{1}{2}} e^{-i t \left( \frac{v^2}{4} + v \xi \right)} \left[ t v \, \widehat{g}\left( t^{\frac{1}{2}} \left( \xi + \frac{v}{2} \right) \right) + i t^{\frac{1}{2}}\p_\xi \widehat{g}\left( t^{\frac{1}{2}} \left( \xi + \frac{v}{2} \right) \right) \right].
\end{align*}
 By Plancherel's Theorem (Theorem \ref{plancherel}):
\[
\| P^+ (x \overline{\Phi_v^+}) \|_{L^2}^2 =  \int_0^\infty \left| \widehat{x \overline{\Phi_v^+}}(\xi) \right|^2 d\xi.
\]
Then,
\[
\left| \widehat{x \overline{\Phi_v^+}}(\xi) \right|^2 \le 2 t \left[ (tv)^2 \left| \widehat{g}\left( t^{\frac{1}{2}} \left( \xi + \frac{v}{2} \right) \right) \right|^2 + t \left| \p_\xi \widehat{g}\left( t^{\frac{1}{2}} \left( \xi + \frac{v}{2} \right) \right) \right|^2 \right].
\]
Next we perform the change of variables $\eta = t^{\frac{1}{2}} \left( \xi + \frac{v}{2} \right)$, so $d\xi = t^{-\frac{1}{2}} d\eta$.  Then for $v>0,$ we have 
\begin{align*}
\| P^+ (x \overline{\Phi_v^+}) \|_{L^2}^2 &\lesssim t \cdot t^{-\frac{1}{2}} \int_{\frac{t^{\frac{1}{2}}v}{2}}^\infty \left[ (tv)^2 |\widehat{g}(\eta)|^2 + t |\p_\xi \widehat{g}(\eta)|^2 \right] d\eta \\
&= \int_{\frac{t^{\frac{1}{2}}v}{2}}^\infty \left[ t^{5/2} v^2 |\widehat{g}(\eta)|^2 + t^{\frac{3}{2}} |\p_\xi \widehat{g}(\eta)|^2 \right] d\eta \\
&= t^{\frac{3}{2}} \left[ (t^{\frac{1}{2}} v)^2 \int_{\frac{t^{\frac{1}{2}}v}{2}}^\infty |\widehat{g}(\eta)|^2 \, d\eta + \int_{\frac{t^{\frac{1}{2}}v}{2}}^\infty |\p_\xi \widehat{g}(\eta)|^2 \, d\eta \right] \\
&\lesssim_N t^{\frac{3}{2}} (t^{\frac{1}{2}} v)^2 \la t^{\frac{1}{2}} v \ra^{-2N} + t^{\frac{3}{2}} \la t^{\frac{1}{2}} v \ra^{-2N} \\
&= t^{\frac{3}{2}} \left( 1 + (t^{\frac{1}{2}} v)^2 \right) \la t^{\frac{1}{2}} v \ra^{-2N} \\
&\lesssim t^{\frac{3}{2}} \la t^{\frac{1}{2}} v \ra^{2(1-N)}.
\end{align*}
Taking the square root yields
\[
\| P^+ (x \overline{\Phi_v^+}) \|_{L^2} \lesssim_N t^{\frac{3}{4}} \la t^{\frac{1}{2}} v \ra^{1-N}.
\]
Similarly for $v<0$ we have
\[
\| P^- (x \overline{\Phi_v^+}) \|_{L^2} \lesssim_N t^{\frac{3}{4}} \la t^{\frac{1}{2}} v \ra^{1-N}.
\]
\end{proof}
Now we begin to bound the remainders.

\subsection*{\texorpdfstring{Estimate of the Term $R_1^+$}{Estimate of the Term R1} }

Define the auxiliary function
\[
g_{v,t}(x) := (x-tv)\chi_{v,t}(x) - 2it^{\frac{1}{2}}\chi'_{v,t}(x),
\]
so that the error term takes the form
\[
R_1^+ = \frac{1}{2t}\int_{\mathbb{R}} e^{-i\phi^+}\px\overline{g_{v,t}}(x)\,u^+(x)\,dx.
\]
Note that using \eqref{chi norms}, we get
\begin{align}\label{g norms}
    \|g_{v,t}\|_{L^2}\lesssim t^{\frac{3}{4}}, \,\,\,     \|\px g_{v,t}\|_{L^2}\lesssim t^{\frac{1}{4}},\,\,\, \left\| \lp \frac{x-vt}{t} \rp g_{v,t} \right\|_{L^2}\lesssim t^{1/4}.
\end{align}
Integrating by parts and using \eqref{dw},  we rewrite $R_1^+$ as:
\begin{align*}
    R_1^+ &= -\frac{1}{2t}\int_{\mathbb{R}} \overline{g_{v,t}}(x)\,\px\left( e^{-i\phi^+}u^+\right) dx \\
    &= \frac{i}{4t^2}\int_{\mathbb{R}} \overline{g_{v,t}}(x)\,e^{-i\phi^+}\left(P^+Lu + c\right) dx.
\end{align*}
We split this expression into two parts.
Using Cauchy--Schwarz and \eqref{g norms}, the first term is bounded by
    \[
    \frac{1}{4t^2}\|g_{v,t}\|_{L^2}\|Lu\|_{L^2} \lesssim t^{-\frac{5}{4}}\|Lu\|_{L^2}.
    \]
 Making the change of variables $y=\frac{x-vt}{\sqrt{t}}$, we obtain
    \[
    \int_{\mathbb{R}} \overline{g_{v,t}}(x)\,e^{-i\phi^+}dx = t\,e^{-i\frac{v^2t}{4}}\int_{\mathbb{R}} e^{-i\frac{v t^{\frac{1}{2}}}{2}y}\,h(y)\,dy = \sqrt{2\pi}\,t\,e^{-i\frac{v^2t}{4}}\,\wh h\left(\frac{vt^{\frac{1}{2}}}{2}\right),
    \]
    where
    \[
    h(y) := e^{-i\frac{y^2}{4}}\left( y\chi(y) + 2i\chi'(y) \right).
    \]
Since $\chi$ is real and even, $h$ is an odd Schwartz function. Thus $\wh h(0)=0$, and for any $N \ge 0$, its Fourier transform satisfies $|\wh h(\eta)| \lesssim_N |\eta|\la\eta\ra^{-N}$. Evaluating at $\eta = \frac{vt^{\frac{1}{2}}}{2}$ yields
\[ \left|\wh h\left(\frac{vt^{\frac{1}{2}}}{2}\right)\right| \lesssim_N |v|t^{\frac{1}{2}}\la t^{\frac{1}{2}}v\ra^{-N}. \]
Combining both parts gives the bound for $|R_1^+|$:
\begin{equation}\label{eq:R1_bound}
    |R_1^+| \lesssim t^{-\frac{5}{4}}\|Lu\|_{L^2} + |c|\,|v|\,t^{-\frac{1}{2}}\la t^{\frac{1}{2}}v\ra^{-N}.
\end{equation}

The second term in \eqref{eq:R1_bound} is integrable in time uniformly with respect to $v$. By the change of variables $\sigma = |v|s^{\frac{1}{2}}$, we have:
\begin{align*}
    \int_1^\infty |v|\,s^{-\frac{1}{2}}\la s^{\frac{1}{2}}v\ra^{-N}ds &\le 2\int_0^\infty\la\sigma\ra^{-N}d\sigma \lesssim 1.
\end{align*}
Also,
\begin{align*}
     \int_1^\infty v^2\,s^{-\frac{1}{2}}\la s^{\frac{1}{2}}v\ra^{-N}ds &\lesssim \min\{|v|,\,|v|^{2-N}\} \lesssim 1,
\end{align*}
which will be useful for $vR_1^+.$
Next, we have
\begin{align*}
    vR_1^+ 
    &= \frac{i}{4t^2}\int_{\mathbb{R}} v\overline{g_{v,t}}(x)\,e^{-i\phi^+}\left(P^+Lu + c\right) dx.
\end{align*}
For $v R_1^+$, we apply \eqref{eq:R1_bound} directly for the second term. For the first term, we decompose the phase factor as
\[
v e^{-i\phi^+} = \frac{x}{t}e^{-i\phi^+} - \frac{x-vt}{t}e^{-i\phi^+}, \quad \text{where} \quad \frac{x}{t}e^{-i\phi^+} = 2i\px(e^{-i\phi^+}).
\]
Then, 
\begin{align*}
    \frac{i}{4t^2}\int_{\mathbb{R}} \overline{g_{v,t}} \, (v e^{-i\phi^+}) \, P^+Lu \, dx =& \underbrace{-\frac{1}{2t^2}\int_{\mathbb{R}} \overline{g_{v,t}} \, \px(e^{-i\phi^+}) \, P^+Lu \, dx}_{I_1} \\
    &\underbrace{-\frac{i}{4t^3}\int_{\mathbb{R}} (x-vt)\overline{g_{v,t}} \, e^{-i\phi^+} \, P^+Lu \, dx}_{I_2}.
\end{align*}
For $I_1$, integrating by parts yields
\[
I_1 = -\frac{1}{2t^2} \int_{\mathbb{R}} \px(e^{-i\phi^+}) \, \overline{g_{v,t}} \, P^+Lu \, dx = \frac{1}{2t^2} \int_{\mathbb{R}} e^{-i\phi^+} \left( \px\overline{g_{v,t}} \, P^+Lu + \overline{g_{v,t}} \, \px P^+Lu \right) dx.
\]
Using \eqref{g norms}, an application of Cauchy--Schwarz yields
\[
|I_1| \lesssim t^{-2}\left( t^{\frac{1}{4}}\|Lu\|_{L^2} + t^{\frac{3}{4}}\|\px Lu\|_{L^2} \right).
\]
For $I_2$, we use the third bound in \eqref{g norms}.
\[
|I_2| \lesssim t^{-2} \cdot t^{\frac{1}{4}} \|Lu\|_{L^2} = t^{-\frac{7}{4}}\|Lu\|_{L^2}.
\]
Combining $|I_1|$ and $|I_2|$ with the constant term bound, we obtain
\begin{align}\label{eq:vR1_bound}
    |vR_1^+| &\lesssim t^{-2}\left( t^{\frac{1}{4}}\|Lu\|_{L^2} + t^{\frac{3}{4}}\|\px Lu\|_{L^2} \right) + t^{-\frac{7}{4}}\|Lu\|_{L^2} + |c|v^2t^{-\frac{1}{2}}\la t^{\frac{1}{2}}v\ra^{-N}.
\end{align}
For $R_2^+$, we have
\begin{align*}
R_2^+=   -t^{-\frac{1}{2}}\int \widebar{P^+\Psi_v^+}  (|u|^2 u)\, dx,
\end{align*}
where $\Psi_v^+:=e^{i\phi^+}\chi'\lp \frac{x-vt}{\sqrt{t}}\rp.$ 
Then,
\begin{align*}
    |R_2^+|
    \lesssim& t^{-\frac{1}{2}}\|u\|_{L^\infty}^2 \|u\|_{L^2}\left\|  \Psi_v^+\right\|_{L^2}\\
\lesssim& t^{-\frac{1}{4}}\|u\|_{L^\infty}^2 \|u\|_{L^2}\\
\leq& D^2\ep^3 t^{-\frac{5}{4}},
\end{align*}
by applying the bootstrap bounds.
Next, we have
\begin{align*}
    vR_2^+=\int vu|u|^2t^{-\frac{1}{2}}\widebar{P^+\Psi_v^+}.
\end{align*}
We will split this term into regions.

\noindent \textbf{Region I: $|v| \leq t^{-\frac12+\delta}$.} In this case, we use the bound on $R_2^+$ to get

\begin{align*}
 |vR_2^+| \leq  |v| |R_2^+|
    \lesssim& D^2\ep^3 t^{-\frac{7}{4}+\delta}.
\end{align*}
Regions II and III are similar, but Region II is more complicated, so we start with Region III.

\noindent\textbf{Region III: $v \leq -t^{-\frac12+\delta}$. } 

In this case, $v<0$ so we can apply \eqref{P- L^2 bound} of Lemma \ref{Packet tails lemma} to get
\begin{align*}
     |vR_2^+|=&\left|\int vu|u|^2t^{-\frac{1}{2}}\widebar{P^+\Psi_v^+} \right|\lesssim v\|u\|_{L^\infty}^2\|u\|_{L^2} t^{-\frac{1}{2}}\|P^-\widebar{\Psi_v^+} \|_{L^2}\\
     \lesssim& D^2\ep^3 t^{-\frac{3}{2}}vt^{\frac{1}{4}} \la t^{\frac{1}{2}}v  \ra^{-N}= D^2\ep^3 t^{-\frac{7}{4}} \la t^{\frac{1}{2}}v  \ra^{1-N}
\end{align*}
 Since $|v|\geq t^{-\frac{1}{2}+\delta},$ we use $(|v|t^{\frac{1}{2}})^{-1}\leq t^{-\delta}$. Taking $N=2,$ we obtain the bound
 \[|vR_2^+|\lesssim D^2\ep^3 t^{-\frac{7}{4}-\delta}.\]

\noindent\textbf{Region II: $v \geq t^{-\frac12+\delta}$.} 

Here we first split up the term, and then bound each piece.
\begin{align*}
     vR_2^+=\int vu|u|^2t^{-\frac{1}{2}}P^-\widebar{\Psi_v^+} =\underbrace{\int vu|u|^2t^{-\frac{1}{2}}\widebar{\Psi_v^+}}_{I} -\underbrace{\int vu|u|^2t^{-\frac{1}{2}}P^+\widebar{\Psi_v^+} }_{II}.
\end{align*}
For the second term, we use \eqref{P+ L^2 bound} of Lemma \ref{Packet tails lemma}.
\begin{align*}
    |II|\lesssim& v\|u\|_{L^\infty}^2\|u\|_{L^2} t^{-\frac{1}{2}}\|P^+\widebar{\Psi_v^+} \|_{L^2}\\
     \lesssim& D^2\ep^3 t^{-\frac{3}{2}}vt^{\frac{1}{4}} \la t^{\frac{1}{2}}v  \ra^{-N}= D^2\ep^3 t^{-\frac{7}{4}} \la t^{\frac{1}{2}}v  \ra^{1-N}.
\end{align*}
Taking $N=2$ as above, we have $|vR_2^+|\lesssim D^2\ep^3 t^{-\frac{7}{4}-\delta}.$  For the first term, we split it as
\begin{align*}
  I=\int vu|u|^2t^{-\frac{1}{2}}\widebar{\Psi_v^+} =\underbrace{\int \frac{x}{t} u(|u|^2)t^{-\frac{1}{2}}\widebar{\Psi_v^+}}_{I_a} +\underbrace{ \int \lp v- \frac{x}{t}\rp u(|u|^2)t^{-\frac{1}{2}}\widebar{\Psi_v^+} }_{I_b}
\end{align*}
where 
\begin{align*}
    |I_b|\lesssim t^{-\frac{1}{2}}\|u\|_{L^\infty}^3 \int \left| \lp v-\frac{x}{t} \rp \widebar{ \Psi^+_v}\right|\approx t^{-\frac{1}{2}}\|u\|_{L^\infty}^3\lesssim D^3\ep^3 t^{-2}.
\end{align*}
Then, 
\begin{align*}
I_a=\int \frac{x}{t} u(|u|^2)t^{-\frac{1}{2}}\widebar{\Psi_v^+} =\int \left(\frac{Lu}{t}+2Hu_x\right) |u|^2t^{-\frac{1}{2}}\widebar{\Psi_v^+}.
\end{align*}
Then the first term is
\begin{align*}
    \left|\int t^{-\frac{3}{2}}Lu|u|^2\widebar{\Psi_v^+} \, dx\right|\lesssim t^{-\frac{3}{2}}\|u\|_{L^\infty}^2\|Lu\|_{L^2}\|\widebar{\Psi_v^+} \|_{L^2},
\end{align*}
where $\|\widebar{\Psi_v^+} \|_{L^2}\leq\|\Psi_v^+\|_{L^2}\approx t^\frac{1}{4}$. The second term is 
\begin{align*}
    \left|\int Hu_x|u|^2t^{-\frac{1}{2}}\widebar{P^+\Psi_v^+}\, dx\right|\lesssim& t^{-\frac{1}{2}}\|H u_x\|_{L^\infty}\|u\|_{L^\infty}\|u\|_{L^2}\|P^+\Psi_v^+\|_{L^2}\\
    \lesssim& t^{-\frac{1}{4}}(\|u_x^+\|_{L^\infty}+\|u_x^-\|_{L^\infty})\|u\|_{L^\infty}\|u\|_{L^2}.
\end{align*}
Overall, we get the bound
\begin{align*}
    |vR_2^+|\lesssim D^2\ep^2t^{-\frac{9}{4}}\|Lu\|_{L^2}+ D^2\ep^3t^{-\frac{5}{4}}
\end{align*}

Next, we bound the term $R_3^+.$  Recall $R_3^+=R_{3,1}^++R_{3,2}^++R_{3,3}^++R_{3,4}^++R_{3,5}^+$, where
$A_i=R_{3,i}^++B_i$ for $i=1,4$ are the resonant terms and $A_j=R_{3,j}^+$ for $j=2,3,5$ are the nonresonant terms.   
For the latter case, we rely on the following Lemma.
\begin{lemma}\label{nonres lemma}
Let $t > 0$, and let $G \in W^{1,1}(\mathbb{R})$ be supported in $\{x \in \mathbb{R} : |x - vt| \leq C_\chi t^{\frac{1}{2}}\}$. Then
\begin{equation*}
    \left| \int_{\mathbb{R}} e^{\pm 2i\phi^+} \frac{x}{t} G(x) \, dx \right| \leq \int_{\mathbb{R}} |\partial_x G(x)| \, dx.
\end{equation*}
\end{lemma}

\begin{proof}
Since $2\phi^+ = \frac{x^2}{2t}$, differentiating yields
\begin{equation*}
    \partial_x \left( e^{\pm 2i\phi^+} \right) = \pm i \frac{x}{t} e^{\pm 2i\phi^+}.
\end{equation*}
Substituting this identity into the integral and integrating by parts gives
\begin{align*}
    \int_{\mathbb{R}} e^{\pm 2i\phi^+} \frac{x}{t} G(x) \, dx 
    &= \mp i \int_{\mathbb{R}} \partial_x \left( e^{\pm 2i\phi^+} \right) G(x) \, dx \\
    &= \pm i \int_{\mathbb{R}} e^{\pm 2i\phi^+} \partial_x G(x) \, dx,
\end{align*}
Therefore
\begin{align*}
    \left| \int_{\mathbb{R}} e^{\pm 2i\phi^+} \frac{x}{t} G(x) \, dx \right|
    &\leq \int_{\mathbb{R}} |\partial_x G(x)| \, dx.
\end{align*}
This completes the proof.
\end{proof}

We now start to bound $R_3^+$  by splitting into regions. Recall that $\delta<\frac{1}{4}$ is fixed.

\noindent \textbf{Region I: $|v| \leq t^{-\frac12+\delta}$.} In this region where $|v|$ is small, we bound $\sum_{j=1}^5 A_j$ directly. Recall 
\begin{align*}
    A_j:=\frac{i}{2t}\int x\,P^+(f_j)\,\widebar{\Phi_v^+}\,dx.
\end{align*}

On $\text{supp}\,\Phi_v^+$, we have $|x| \leq C_\chi \sqrt{t} + |v|t$. Then,
\begin{align*}
    |A_j| &\leq \frac{1}{2t} \| P^+ f_j\|_{L^2} \|x \Phi_v^+\|_{L^2} \\
    &\lesssim t^{-1} \|u\|_{L^\infty}^2 \|u\|_{L^2} \left( C_\chi \sqrt{t} + |v|t \right) t^{1/4} \\
    &\lesssim D^2 \ep^3 t^{-2} \left( C_\chi \sqrt{t} + |v|t \right) t^{1/4} \\
    &\lesssim D^2 \ep^3 \left( C_\chi t^{-\frac{5}{4}} + t^{-\frac{5}{4}+\delta} \right) \lesssim D^2 \ep^3 t^{-\frac54+\delta}.
\end{align*}
This finishes the bound for $j=2,3,5.$  For $j=1,4,$ we must also bound $B_1, B_4$ to get the bound for $R_{3,1}^+, R_{3,4}^+.$ Note that applying the bootstrap bound to \eqref{gamma+bnds} gives $|\gamma^+| \lesssim D \ep$. Using this along with $|v| \leq t^{-\frac{1}{2}+\delta}$ gives
\begin{align*}
    |B_1+B_4| &\leq \frac{3|v|}{2} \|u\|_{L^\infty}^2 |\gamma^+| \\
    &\lesssim t^{-\frac{1}{2}+\delta} \left( D \ep t^{-\frac{1}{2}} \right)^2 (D \ep) \\
    &\lesssim D^3 \ep^3 t^{-\frac{3}{2}+\delta}.
\end{align*}
Thus, in Region I, we have
\begin{align*}
    |R_3^+|\lesssim D^2 \ep^3 t^{-\frac54+\delta}+ D^3 \ep^3 t^{-\frac{3}{2}+\delta}.
\end{align*}
Also,
\begin{align*}
    |vR_3^+|\lesssim D^2 \ep^3 t^{-\frac74+2\delta}+ D^3 \ep^3 t^{-2+2\delta}.
\end{align*}

\noindent\textbf{Region II:} $v \geq t^{-\frac12+\delta}$ (Nonresonant terms $j \in \{2,3,5\}$).

We use the adjoint property to write
\begin{align*}
        A_j=\frac{i}{2t}\int x\,P^+(f_j)\,\widebar{\Phi_v^+}\,dx=    \frac{i}{2t}\int f_j\,P^-(x\widebar{\Phi_v^+})\,dx.
\end{align*}
 We then split $A_j$ using $P^-(x\widebar{\Phi_v^+})=x\widebar{\Phi_v^+}-P^+(x\widebar{\Phi_v^+})$, into a main localized integral $\tilde{A}_j$ and a projection tail $E_j$:
\begin{equation*}
    A_j = \tilde{A}_j + E_j, \quad \text{where} \quad \tilde{A}_j := \frac{i}{2t} \int_{\mathbb{R}} x f_j \overline{\Phi_v^+} \, dx, \quad E_j := -\frac{i}{2t} \int_{\mathbb{R}} f_j P^+\left( x \overline{\Phi_v^+} \right) dx.
\end{equation*}
By Lemma \ref{Packet tails lemma}, since $v\geq t^{-\frac12+\delta},$ choosing $N=2$ gives
\begin{align*}
    |E_j|\lesssim t^{-1}\|u\|_{L^\infty}^2\|u\|_{L^2}\|P^+\left( x \overline{\Phi_v^+} \right) \|_{L^2}\leq D^2\ep^3t^{-2}\cdot t^{\frac{3}{4}}\la t^{\frac{1}{2}}v
\ra^{1-N}
\lesssim D^2\ep^3 t^{-\frac{5}{4}-\delta}.
\end{align*}
Also, for $N=3,$
\begin{align*}
    |vE_j|\lesssim vt^{-1}\|u\|_{L^\infty}^2\|u\|_{L^2}\|P^+\left( x \overline{\Phi_v^+} \right) \|_{L^2}\leq D^2\ep^3t^{-2}\cdot t^{\frac{3}{4}}v\la t^{\frac{1}{2}}v \ra^{1-N}\lesssim D^2\ep^3t^{-\frac{7}{4}-\delta}.
\end{align*}
For the main terms $\tilde{A}_j$, we demonstrate the argument for $j = 2$, and the cases $j = 3, 5$ are identical. Writing $f_2\widebar{\Phi_v^+} = (u^+)^2 \widebar{\um} \widebar{\Phi_v^+} = e^{2i\phi^+} (w^+)^2 \widebar{w^-} \chi_{v,t}$, we have
\begin{equation*}
    \tilde{A}_2 = \frac{i}{2} \int_{\mathbb{R}} e^{2i\phi^+} \frac{x}{t} G(x) \, dx, \quad \text{where} \quad G(x) := (w^+)^2 \widebar{w^-} \chi_{v,t}.
\end{equation*}
Applying Lemma \ref{nonres lemma} yields
\begin{equation*}
    |\tilde{A}_2| \leq \frac{1}{2} \int_{\mathbb{R}} |\partial_x G(x)| \, dx.
\end{equation*}
Expanding $\partial_x G$ via the product rule gives
\begin{align*}
    \int_{\mathbb{R}} |\partial_x G| \, dx 
    &\leq 2 \|w^+\|_{L^\infty} \|\bar{w}^-\|_{L^\infty} \int_{\mathbb{R}} |\partial_x w^+| \chi_{v,t} \, dx 
    + \|w^+\|_{L^\infty}^2 \int_{\mathbb{R}} |\partial_x \bar{w}^-| \chi_{v,t} \, dx \\
    &\quad + \|w^+\|_{L^\infty}^2 \|\bar{w}^-\|_{L^\infty} \|\partial_x \chi_{v,t}\|_{L^1}.
\end{align*}
Next, recall that
\begin{align*}
    \|w^\pm\|_{L^\infty} \lesssim& \|u\|_{L^\infty}\lesssim D\ep t^{-\frac{1}{2}}.
\end{align*}
Also, since $\px \wp =-\frac{i}{2t}e^{-i\phi^+}(P^+ Lu +c),$ we have 
\begin{align*}
    \int_{\mathbb{R}} |\partial_x w^+| \chi_{v,t} \, dx \leq t^{-1}\lp  \|Lu\|_{L^2} \|\chi_{v,t}\|_{L^2}+ |c| \|\chi_{v,t}\|_{L^1}\rp \leq t^{-1}\lp t^{1/4} \|Lu\|_{L^2} + t^{\frac{1}{2}} |c| \rp,
\end{align*}
and the same for $\wm.$  Putting this all together, we obtain
\begin{equation}\label{eq:nonres_bound}
    |\tilde{A}_j| \lesssim D^2 \ep^2 \lp t^{-\frac74} \|Lu\|_{L^2}+t^{-\frac{3}{2}}|c| \rp+ D^3 \ep^3 t^{-\frac32}, \qquad j \in \{2, 3, 5\}.
\end{equation}

To bound $v \tilde{A}_2$, we decompose the velocity weight as $v = \frac{x}{t} - \frac{x - vt}{t}$:
\begin{align*}
v \tilde{A}_2&= \frac{i}{2} \int e^{2i\phi^+} \frac{x}{t} \cdot \frac{x}{t} G(x) \, dx + \frac{i}{2} \int e^{2i\phi^+} \frac{x}{t} \left( \frac{vt-x}{t} \right) G(x) \, dx \\[10pt]
&= \frac{1}{2} \int \partial_x \left(e^{2i\phi^+}\right) G(x) \, dx + \frac{1}{2} \int \partial_x \left(e^{2i\phi^+}\right) \cdot \left( \frac{vt-x}{t} \right) G(x) \, dx \\[10pt]
&= -\frac{1}{2} \int e^{2i\phi^+} \partial_x G(x) \, dx - \frac{1}{2} \int e^{2i\phi^+} \partial_x \left[ \left(\frac{vt-x}{t}\right) G(x) \right] dx.
\end{align*}
Then, the first term can be bounded as for $\tilde{A}_2$ and for the last term we use that $|x - vt| \le C_\chi t^{\frac{1}{2}}$ on $\supp\chi_{v,t}$:
\begin{align*}
|v \tilde{A}_2|
\lesssim&  \int |\partial_x G(x)| \, dx + \; \frac{1}{t} \int |G(x)| \, dx + \int \left| \left( \frac{vt-x}{t} \right) \partial_x  G(x)  \right| dx \\[10pt]
\lesssim&  D^2 \ep^2 \lp t^{-\frac74} \|Lu\|_{L^2}+t^{-\frac{3}{2}}|c| \rp+ D^3 \ep^3 t^{-\frac32} + t^{-1} \|u\|_{L^\infty}^3 \|\chi_{v,t}\|_{L^1}\\
&+ \; C_\chi t^{-1/2} \lp D^2 \ep^2 \lp t^{-\frac74} \|Lu\|_{L^2}+t^{-\frac{3}{2}}|c| \rp+ D^3 \ep^3 t^{-\frac32} \rp.
\end{align*}
Condensing this, we obtain the bound
\begin{align*}
    |v \tilde{A}_j|\lesssim D^2 \ep^2 \lp t^{-\frac74} \|Lu\|_{L^2}+t^{-\frac{3}{2}}|c| \rp+ D^3 \ep^3 t^{-\frac32}.
\end{align*}
for $j \in \{2,3,5\}$.

\noindent\textbf{Region II:} $v \geq t^{-\frac12+\delta}$ (Resonant terms $j \in \{1,4\}$).

Next we need to bound $R_{3,1}^+$ and $R_{3,4}^+$ . Recall that for $i=1,4$ $A_i = R_{3,i}^+ + B_i.$
Just as we did for the nonresonant terms, in this case we can also write:
\begin{equation*}
A_i = \widetilde{A}_i + E_i,  
\end{equation*}
where $\widetilde{A}_i, E_i$ are as defined above. Then, $R_{3,i}^+ = \widetilde{A}_i + E_i - B_i$ for $i=1,4$.
We again have
\begin{equation*}
\|E_{i}\|_{L^\infty} \lesssim D^2 \epsilon^3 t^{-\frac{5}{4} - \delta}, \qquad \|v E_{i}\|_{L^{\infty}} \lesssim D^2 \epsilon^3 t^{-\frac{7}{4} - \delta}.
\end{equation*}
Then consider $\widetilde{A}_i - B_i$. Without loss of generality, let $i=1$, since the proof for $i=4$ is identical. Then,
\begin{align*}
\widetilde{A}_1 - B_1
&= \frac{i}{2} \int \frac{x}{t} |u^+|^2 u^+ \overline{\Phi_v^+} \, dx - \frac{i v}{2} |u^+(t,vt)|^2 \gamma^+(t,v) \\
&= \frac{i}{2} \int \frac{x}{t} |u^+|^2 u^+ \overline{\Phi_v^+} \, dx - \frac{i v}{2} |u^+(t,vt)|^2 \int u^+(t,x) \overline{\Phi_v^+} \, dx \\
&= \frac{i}{2} \int \left( \frac{x}{t} |w^+(t,x)|^2 - v |w^+(t,vt)|^2 \right) {w^+(t,x)} \chi_{v,t} \, dx.
\end{align*}
By adding and subtracting $\frac{i}{2} \int v |w^+(t,x)|^2 w^+(t,x) \chi_{v,t} \, dx$, we obtain
\begin{align*}
\widetilde{A}_1 - B_1&= \frac{i}{2} \int \left( \frac{x}{t} - v \right) |w^+(t,x)|^2 w^+(t,x) \chi_{v,t} \, dx \\
&\quad + \frac{i}{2} \int v \left( |w^+(t,x)|^2 - |w^+(t,vt)|^2 \right) w^+(t,x) \chi_{v,t} \, dx.
\end{align*}

\noindent Since $|x - vt| \le C_\chi t^{\frac{1}{2}}$ on $\supp \chi_{v,t}$, the first integral is bounded by

\begin{equation*}
C t^{-\frac{1}{2}} \|u\|_{L^\infty}^3 \|\chi_{v,t}\|_{L^1} \lesssim (D\epsilon)^3 t^{-\frac{3}{2}}.
\end{equation*}
For the second, we apply Lemma \ref{local diff lemma} to get the bound
\begin{align*}
   & \|u\|_{L^\infty}^2 \sup_v\||v| (\wp(t,x)-\wp(t,vt))\|_{L^\infty}\|\chi_{v,t}\|_{L^1} \lesssim  D^2\ep^2 t^{-\frac{5}{4}}(\|Lu\|_{H^1}+|c|).
\end{align*}
Overall we have
\begin{align*}
    |\widetilde{A}_1 - B_1|+|\widetilde{A}_4 - B_4|\lesssim D^3\epsilon^3 t^{-\frac{3}{2}}+ D^2\ep^2 t^{-\frac{5}{4}}(\|Lu\|_{H^1}+|c|).
\end{align*}
Next,
\begin{align*}
v\widetilde{A}_1 - vB_1=& \frac{i}{2} \int v\left( \frac{x}{t} - v \right) |w^+(t,x)|^2 w^+(t,x) \chi_{v,t} \, dx \\
&+ \frac{i}{2} \int v^2 \left( |w^+(t,x)|^2 - |w^+(t,vt)|^2 \right) w^+(t,x) \chi_{v,t} \, dx\\
=&\frac{i}{2} \int \frac{x}{t}\left( \frac{x}{t} - v \right) |w^+(t,x)|^2 w^+(t,x) \chi_{v,t} \, dx \\
&-\frac{i}{2} \int \left( \frac{x}{t} - v \right)^2 |w^+(t,x)|^2 w^+(t,x) \chi_{v,t} \, dx \\
&+ \frac{i}{2} \int v\cdot\frac{x}{t} \left( |w^+(t,x)|^2 - |w^+(t,vt)|^2 \right) w^+(t,x) \chi_{v,t} \, dx\\
&- \frac{i}{2} \int v\left( \frac{x}{t} - v \right)  \left( |w^+(t,x)|^2 - |w^+(t,vt)|^2 \right) w^+(t,x) \chi_{v,t} \, dx\\
=&I+II+III+IV.
\end{align*}
 For the first term, we combine $\px( \wp)=\frac{-ix}{2t}e^{-i\phi^+} \up +e^{-i\phi^+}  \px \up$ with \eqref{dw} to get the identity
\begin{align} \label{x/t w identity}
    \frac{x}{t} w^+ = \frac{1}{t} e^{-i\phi^+} (P^+ Lu + c) - 2 i j^+.
\end{align}
Then the first term becomes
\begin{align*}
    I=\frac{i}{2} \int \left( \frac{x}{t} - v \right) |w^+(t,x)|^2 \lp\frac{1}{t} e^{-i\phi^+} (P^+ Lu + c) - 2 i j^+\rp \chi_{v,t}(x) \, dx,
\end{align*}
which is bounded as
\begin{align*}
    |I|\lesssim& t^{-\frac{1}{2}} \|u\|_{L^\infty}^2 \left[\int \frac{1}{t}(|P^+ Lu| + |c|)\chi_{v,t} \, dx+\int  | j^+|)\chi_{v,t} \, dx\right]\\
   \lesssim& D^2\ep^2 t^{-\frac{3}{2}} \lp t^{-1}\|Lu\|_{L^2} \|\chi_{v,t}\|_{L^2} +t^{-1}|c|\|\chi_{v,t}\|_{L^1}+\|u_x\|_{L^\infty}\|\chi_{v,t}\|_{L^1}  \rp \\
   \lesssim& D^2\ep^2 t^{-\frac{3}{2}} \lp t^{-\frac{3}{4}}\|Lu\|_{L^2}  +t^{-\frac{1}{2}}|c|+D\epsilon \rp.
\end{align*}
Also
\begin{align*}
    |II|\lesssim t^{-2} D^3\ep^3.
\end{align*}
Next, using Lemma \ref{local diff lemma} and  \eqref{x/t w identity}, we have
\begin{align*}
    |III|\lesssim& \|u\|_{L^\infty}\sup_v\||v| (\wp(t,x)-\wp(t,vt))\|_{L^\infty}\int \lp\frac{1}{t} e^{-i\phi^+} (P^+ Lu + c) - 2 i j^+\rp \chi_{v,t} \, dx \\
 \lesssim& D\ep  t^{-\frac54}\lp\|Lu\|_{H^1}+|c|\rp \lp t^{-\frac{3}{4}}\|Lu\|_{L^2}  +t^{-\frac{1}{2}}|c|+D\epsilon \rp.
\end{align*}
Also,
\begin{align*}
    |IV|\lesssim& D^2\ep^2 t^{-1} \sup_v\||v| (\wp(t,x)-\wp(t,vt))\|_{L^\infty}\\
    \lesssim&D^2\ep^2 t^{-\frac{7}{4}} \lp\|Lu\|_{H^1}+|c|\rp.
\end{align*}
Overall, we have
\begin{align*}
     |v\widetilde{A}_1 - vB_1|+|v\widetilde{A}_4 -v B_4|
     \lesssim& D^3\ep^3 t^{-\frac{3}{2}}+D^2\ep^2  t^{-\frac54}\lp\|Lu\|_{H^1}+|c|\rp+D\ep t^{-\frac74} \lp\|Lu\|_{H^1}+|c|\rp^2.
\end{align*}
Putting this all together, in Region II we have
\begin{align*}
|R_3^+|\lesssim&  D^3\epsilon^3 t^{-\frac{3}{2}}+ D^2\ep^2 t^{-\frac{5}{4}}(\|Lu\|_{H^1}+|c|).\\
 |vR_3^+|\lesssim&   D^2\ep^2  t^{-\frac54}\lp\|Lu\|_{H^1}+|c|\rp+D^2 \ep^2 t^{-\frac32}( |c| + D\ep)+D\ep t^{-\frac74} \lp\|Lu\|_{H^1}+|c|\rp^2.
\end{align*}

\noindent\textbf{Region III:} $v \leq -t^{-\frac12+\delta}$.

\noindent For region III, everything is a tail so it is much easier than the previous case.  Consider
\begin{align*}
        A_j=\frac{i}{2t}\int x\,P^+(f_j)\,\widebar{\Phi_v^+}\,dx=    \frac{i}{2t}\int f_j\,P^-(x\widebar{\Phi_v^+})\,dx.
\end{align*}
Using Lemma \ref{Packet tails lemma}, since $v\leq -t^{-\frac12+\delta},$ by choosing $N=2$ we have
\begin{align*}
    |A_j|\lesssim t^{-1}\|u\|_{L^\infty}^2\|u\|_{L^2}\|P^-\left( x \overline{\Phi_v^+} \right) \|_{L^2}\leq D^2\ep^3t^{-2}\cdot t^{\frac{3}{4}}\la t^{\frac{1}{2}}v
\ra^{1-N}
\lesssim D^2\ep^3 t^{-\frac{5}{4}-\delta}.
\end{align*}
Also, for $N=3,$
\begin{align*}
    |vA_j|\lesssim vt^{-1}\|u\|_{L^\infty}^2\|u\|_{L^2}\|P^-\left( x \overline{\Phi_v^+} \right) \|_{L^2}\leq D^2\ep^3t^{-2}\cdot t^{\frac{3}{4}}v\la t^{\frac{1}{2}}v \ra^{1-N}\lesssim D^2\ep^3t^{-\frac{7}{4}-\delta}.
\end{align*}
This finishes the bounds for $R_{3,j}^+$ for $j=2,3,5.$ We now move to the resonant terms. 
In this region, we have $|v|\geq t^{-\frac{1}{2}+\delta}$, so consider
\begin{align*}
   | B_1|=&\left| v|\up|^2 \gp\right|=\left| v|\up|^2 \int \up(t,x) \widebar{\Phi_v^+}\, dx\right| \lesssim |v|\|u\|_{L^\infty}^2 \|u\|_{L^2} \|P^- \widebar{\Phi_v^+}\|_{L^2}\\
    \lesssim& D^2\ep^3t^{-1} |v|t^{1/4}\la t^{\frac{1}{2}}v \ra^{-N}\leq D^2\ep^3t^{-\frac{5}{4}} \la t^{\frac{1}{2}}v \ra^{1-N} \leq D^2\ep^3t^{-\frac{5}{4}-\delta}.
\end{align*}
Also
\begin{align*}
   |v B_1|
    \lesssim& D^2\ep^3t^{-1} |v|^2t^{1/4}\la t^{\frac{1}{2}}v \ra^{-N}\leq D^2\ep^3t^{-\frac{7}{4}} \la t^{\frac{1}{2}}v \ra^{2-N} \leq D^2\ep^3t^{-\frac{7}{4}-\delta}.
\end{align*}
Similarly, we have $|B_4|\lesssim D^2\ep^3t^{-\frac{5}{4}-\delta}, |vB_4|\lesssim  D^2\ep^3t^{-\frac{7}{4}-\delta}$. 
Altogether, in Region III, we have the bounds
\begin{align*}
    |R_3^+|\lesssim& D^2\ep^3 t^{-\frac{5}{4}-\delta}\\
  |vR_3^+|\lesssim& D^2\ep^3t^{-\frac{7}{4}-\delta}.
\end{align*}
Combining all the regions, we overall have the bounds:
\begin{align*}
     |R_3^+|\lesssim D^2 \ep^3 t^{-\frac54+\delta}+ D^3 \ep^3 t^{-\frac{3}{2}+\delta}+D^2\ep^2 t^{-\frac{5}{4}}(\|Lu\|_{H^1}+|c|)
\end{align*}
\begin{align*}
    |vR_3^+|\lesssim D^2 \ep^3 t^{-\frac74+2\delta}+ D^3 \ep^3 t^{-2+2\delta}+ D^2\ep^2  t^{-\frac54}\lp\|Lu\|_{H^1}+|c|\rp+D\ep t^{-\frac74} \lp\|Lu\|_{H^1}+|c|\rp^2.
\end{align*}
Lastly, we bound $R_4^+.$  Recall
\begin{align*}
   R_4^+=& i\frac{v}{2}\Big[\lp |\up(t,vt)|^2-t^{-1}|\gp|^2\rp+2\lp |\um(t,vt)|^2-t^{-1}|\gm|^2 \rp \Big]\gp(t,v).
\end{align*}
Since $|u^\pm(t,vt)|=|w^\pm(t,vt)|$ and $t^{-\frac12}\int\chi_{v,t}\,dx=1$, we have
\begin{align*}
    t^{-\frac12}\gamma^\pm-w^\pm(vt)=t^{-\frac12}\int(w^\pm(x)-w^\pm(vt))\chi_{v,t}\,dx.
\end{align*}
This allows us to write each term in $R_4^+$ as
\begin{align*}
  \left|  i\frac{v}{2}\lp |w^\pm(t,vt)|^2-t^{-1}|\gamma^\pm|^2\rp\gp\right|\lesssim& |\gp|  (\|u\|_{L^\infty}+t^{-\frac{1}{2}}\|\gamma\|_{L^\infty}) t^{-\frac12} \left||v|\int(w^\pm(x)-w^\pm(vt))\chi_{v,t}\,dx\right|\\
  \lesssim& D\ep \cdot D\ep t^{-1} \sup_{|x-vt|\leq C_\chi t^{\frac{1}{2}}}|v|\,|w^\pm(x)-w^\pm(vt)|\,\|\chi_{v,t}\|_{L^1}\\
    \lesssim& D^2\ep^2 t^{-\frac{1}{2}}\cdot t^{-\frac{3}{4}}(\|Lu\|_{H^1}+|c|),
\end{align*}
where in the last line we used \eqref{local diff 2} of Lemma \ref{local diff lemma}.
Therefore 
\begin{align*}
    |R_4^+|\lesssim D^2\ep^2  t^{-\frac{5}{4}}(\|Lu\|_{H^1}+|c|).
\end{align*}
Similarly, using  \eqref{v gamma+ L^infty}, we have $\|v\gamma\|_{L^\infty}\lesssim t^{-\frac{3}{4}}\|Lu\|_{L^2_x}+t^{-\frac12}|c|+ D\ep,$ which gives
\begin{align*}
    |vR_4^+|\lesssim D\ep  t^{-\frac{5}{4}}(\|Lu\|_{H^1}+|c|) \lp t^{-\frac{3}{4}}\|Lu\|_{L^2_x}+t^{-\frac12}|c|+ D\ep \rp.
\end{align*}

This finishes the bounds for $R$ and $vR$, which are
\begin{align*}
    |R^+|\lesssim& t^{-\frac{5}{4}}\|Lu\|_{L^2} + |c|\,|v|\,t^{-\frac{1}{2}}\la t^{\frac{1}{2}}v\ra^{-N}+ D^2 \ep^3 t^{-\frac54+\delta}\\
    +& D^3 \ep^3 t^{-\frac{3}{2}+\delta}+D^2\ep^2 t^{-\frac{5}{4}}(\|Lu\|_{H^1}+|c|)\\
    |vR^+| \lesssim&D^3\ep^3 t^{-\frac32}+D^2\ep^3t^{-\frac{5}{4}}+ D^3 \ep^3 t^{-2+2\delta}+  |c|v^2t^{-\frac{1}{2}}\la t^{\frac{1}{2}}v\ra^{-N}\\
    &+  (D^2\ep^2  t^{-\frac54}+t^{-\frac{5}{4}})\lp\|Lu\|_{H^1}+|c|\rp+ D\ep t^{-\frac54} \lp\|Lu\|_{H^1}+|c|\rp^2.
\end{align*}

\subsection{Negative frequency Equation}
To determine the negative frequency equation, we first note that $\gm(u)=\overline{\gp (\bu)}$ since
\begin{align*}
    \overline{\gp (\bu)}=\overline{\int (\bu)^+ \overline{\Phi_v^+}\, dx}=\int \overline{(\bu)^+ }\Phi_v^+\, dx=\int \um \overline{\Phi_v^-}\, dx =\gm(u).
\end{align*}
Then, we have
\begin{align}
    \pt \gp(\bu) =&  i\frac{v}{2t}\lp|\gp(\bu)|^2\gp(\bu)  +2|\gm(\bu)|^2 \gp(\bu)\rp+R^+(\bu).
\end{align}
Taking the conjugate,
\begin{align*}
    \overline{\pt \gp(\bu) }=& \overline{ i\frac{v}{2t}\lp|\gp(\bu)|^2\gp(\bu)  +2|\gm(\bu)|^2 \gp(\bu)\rp+R^+(\bu)}\\
    =& -i\frac{v}{2t}\lp|\overline{\gp(\bu)}|^2\overline{\gp(\bu)}  +2|\overline{\gm(\bu)}|^2 \overline{\gp(\bu)}\rp+\overline{R^+(\bu)}\\
   \pt \gm(u) =& -i\frac{v}{2t} \lp |\gm(u)|^2\gm(u) +2 |\gp(u)|^2\gm(u) \rp+\overline{R^+(\bu)}.
\end{align*}
Define $R^-(u):=\overline{R^+(\bu)}$.  Then the negative frequency equation is 
\begin{align*}
      \pt \gm =& -i\frac{v}{2t} \lp |\gm|^2\gm +2 |\gp|^2\gm \rp+R^-,
\end{align*}
where 
\begin{align*}
    |R^-|\lesssim& t^{-\frac{5}{4}}\|Lu\|_{L^2} + |c|\,|v|\,t^{-\frac{1}{2}}\la t^{\frac{1}{2}}v\ra^{-N}+ D^2 \ep^3 t^{-\frac54+\delta}\\
    +& D^3 \ep^3 t^{-\frac{3}{2}+\delta}+D^2\ep^2 t^{-\frac{5}{4}}(\|Lu\|_{H^1}+|c|)\\
    |vR^-| \lesssim&D^3\ep^3 t^{-\frac32}+D^2\ep^3t^{-\frac{5}{4}}+ D^3 \ep^3 t^{-2+2\delta}+  |c|v^2t^{-\frac{1}{2}}\la t^{\frac{1}{2}}v\ra^{-N}\\
    &+  (D^2\ep^2  t^{-\frac54}+t^{-\frac{5}{4}})\lp\|Lu\|_{H^1}+|c|\rp+ D\ep t^{-\frac54} \lp\|Lu\|_{H^1}+|c|\rp^2.
\end{align*}

\subsection{Asymptotic equation bounds}

Now that we have the remainder bounds, we can find the bounds for $\gamma^\pm.$ We have the equation
\begin{equation}\label{gp eqn again}
    \pt \gp =  i\frac{v}{2t}(|\gp|^2\gp +2|\gm|^2 \gp)+R^+(t,v).
\end{equation}
We perform the energy estimate in the usual way, first by multiplying by $ \widebar{\gp}$,
\begin{align*}
      \pt \gp \cdot \widebar{\gp}=&  i\frac{v}{2t}(|\gp|^4 +2|\gm|^2 |\gp|^2)+R^+(t,v)\cdot \widebar{\gp},
\end{align*}
and then taking the real part of both sides to get
\begin{align} \label{gp EE}
    \frac{d}{dt}|\gp(t,v)|^2=&2\Rea \lp R^+(t,v)\cdot \widebar{\gp} \rp.
\end{align}
To avoid dividing by zero, let $\eta>0.$ Then
\begin{align*}
    \pt \lp |\gp(t)|^2+ \eta^2\rp^{\frac{1}{2}}=\frac{\pt |\gp (t)|^2}{2\lp |\gp(t)|^2+ \eta^2\rp^{\frac{1}{2}}}\leq \frac{ |R^+(t,v)| | \gp| }{\lp |\gp(t)|^2+ \eta^2\rp^{\frac{1}{2}}}\leq |R^+(t,v)|.
\end{align*}
Integrating gives
\begin{align*}
    \lp |\gp(t)|^2+ \eta^2\rp^{\frac{1}{2}} \leq \lp |\gp(1)|^2+ \eta^2\rp^{\frac{1}{2}}+ \int_1^t |R^+(s,v)| \, ds, \text{ for all } \eta>0.
\end{align*}
Sending $\eta \to 0$ gives
\begin{align}\label{gamma+ bound}
     |\gp(t)| \leq |\gp(1)|+ \int_1^t |R^+(s,v)| \, ds,
\end{align}
where $|\gp(1,v)|\leq\|u^+(1)\|_{L^2}\|\Phi_v^+(1)\|_{L^2}\lesssim\ep.$ Next, multiplying by $v$ gives 
\begin{align}\label{vgamma+ bound}
       |v\gp(t)| \lesssim |v\gp(1)|+ \int_1^t |vR^+(s,v)| \, ds.
\end{align}
At $t=1$, \eqref{v gamma+ L^infty}  gives $|v\gp(1,v)|\lesssim\|Lu(1)\|_{L^2}+|c|+\|u^+_x(1)\|_{L^\infty}.$ 
Then, by applying Proposition \ref{vector field estimate general f} to $u_x$, we have
\begin{align*}
  \|u_x^+(1)\|_{L^\infty}^2\leq&\|u_x^+(1)\|_{L^2}\|L(u_x^+)(1)\|_{L^2}.
\end{align*}
By Lemma \ref{commutator lemma}, we have $L(u^\pm_x)=P^\pm L(u_x)$, and so $\|u_x^+(1)\|_{L^\infty}^2\leq\|u_x(1)\|_{L^2}\|L(u_x)(1)\|_{L^2}$.  In the next section, we prove vector field bounds which imply $\|Lu(1)\|_{L^2}\lesssim \ep,\, \|L(u_x)(1)\|_{L^2}\lesssim \ep.$ The local well-posedness of \eqref{mBO} in $H^1$ gives $\|u_x^+(1)\|_{L^2}\lesssim \ep.$
Therefore we get $|v\gp(1,v)|\lesssim \ep.$

\vspace{1em}

The bounds for $R^{\pm}$ depend on the vector field bound, so we turn to that next. 
\section{Vector field equations and bounds \label{Vector field bnds section}}
We now obtain an energy estimate for $L u.$ Recall the definition
\[L=x-2tH\px .\]
First, we find an equation $Lu$ satisfies, and then attempt to find the energy estimate. Applying $L$ to \eqref{mBO}, we have
\begin{align*}
    L( \pt+H\partial_x^2)u =L\partial_{x}(|u|^2u).
\end{align*}
Here, $L$ commutes with the linear part of the equation, and we can use the commutator $[\px, L]=1$, to get 
\begin{align*}
    ( \pt+H\partial_x^2)Lu =\partial_{x}L(|u|^2u)-|u|^2u.
\end{align*}
The problem now is how $L$ acts on $|u|^2u.$  For NLS and DNLS with $J=x+2ti\px,$ we have the property $J(u\bv w)=Ju \bv w-u \widebar{Jv}w+u \bv Jw$.  However, because of the Hilbert transform in \eqref{mBO} and therefore also in our vector field $L,$ we do not have any such property.  Therefore, it is extremely difficult to obtain energy estimates for $Lu$ from this equation.  In order to get the bound, we need to apply a correction to our vector field.

\subsection{\texorpdfstring{Introduce $L^{NL}u$}{Introduce LNLu} }
The idea to use a nonlinear correction to the scaling vector field first appears in \cite{cite:BO} for Benjamin-Ono and is later used in \cite{cite:kdv} with the goal of obtaining dispersive estimates. The concept is explained in a more general setting in the expository work of the Ifrim and Tataru in \cite{cite:WavePax}.

To find the vector field correction, which we will call $L^{NL},$ we first need to go back to how we found $L$ in the first place.  Recall that $Lu$ is the scaling vector field for the linear equation.  If we use this same process, but keep the nonlinear terms included, we obtain the correction vector field, $L^{NL}.$

Recall the scaling symmetry for \eqref{mBO},
\begin{align*}
    u_\lambda(t,x)=\lambda^{\frac{1}{2}}u( \lambda^2 t,\lambda x).
\end{align*}
Let 
\begin{align*}
    v:=\frac{d}{d \lambda}|_{\lambda=1}u_\lambda =\frac{1}{2}u+2tu_t+xu_x.
\end{align*}
Here and in the rest of this section $v$ denotes this function, not a velocity.

\textbf{Claim:} $v$ solves the linearized equation
\[(\pt +H\px^2)v=\px(2|u|^2v+u^2\bv).\]
\begin{proof}
By calculating derivatives of $v,$ we have
\begin{align*}
v_t &= \frac{1}{2} u_t + 2u_t + 2t u_{tt} + x u_{tx} = \frac{5}{2} u_t + 2t u_{tt} + x u_{tx} \\
v_x &= \frac{1}{2} u_x + 2t u_{tx} + u_x + x u_{xx} \\
H v_{xx} &= \frac{3}{2} H u_{xx} + 2t H u_{txx} + H(u_{xx} + x u_{xxx}) \\
&= \frac{5}{2} H u_{xx} + 2t H u_{txx} + H(x u_{xxx}).
\end{align*}
Then, we add them together and use the original equation \eqref{mBO}, resulting in
\begin{align*}
    (\partial_t + H\partial_x^2)v &= \frac{5}{2}(u_t + H u_{xx}) + 2t(u_{tt} + H u_{txx}) + H(x u_{xxx}) + x u_{tx} \\
&= \frac{5}{2} \partial_x(|u|^2 u) + 2t \partial_t(\partial_x(|u|^2 u)) + x \partial_x(\partial_x(|u|^2 u)).
\end{align*}
Here, we were able to commute $H$ with $x$ since $\int u_{xxx}=0.$ To obtain the right hand side of the linearized equation, we compute
\begin{align*}
2|u|^2 v &= |u|^2 u + 4t |u|^2 u_t + 2x |u|^2 u_x \\
u^2 \overline{v} &= \frac{1}{2} |u|^2 u + 2t u^2 \overline{u}_t + x \overline{u}_x u^2.
\end{align*}
From the product rule, we have the commutator
\begin{align*}
[x, \partial_x] f &= x \partial_x f - \partial_x(xf) = -f,
\end{align*}
which we can use to get
\begin{align*}
(\partial_t + H\partial_x^2)v &= \frac{5}{2} \partial_x(|u|^2 u) + \partial_x(2t \partial_t(|u|^2 u)) + \partial_x(x \partial_x(|u|^2 u)) - \partial_x(|u|^2 u) \\
&= \partial_x \left[ \frac{3}{2} |u|^2 u + 2t \partial_t(|u|^2 u) + x \partial_x(|u|^2 u) \right]\\
&=\px \left[  2|u|^2 v+u^2\bv \right].
\end{align*}
This proves the claim. 
\end{proof}
Next, we expand $v$ using the original equation. 
\begin{align*}
v &= \frac{1}{2} u + 2t u_t + x u_x \\
&= \frac{1}{2} u + 2t \left[ \partial_x(|u|^2 u) - H\partial_x^2 u \right] + x u_x\\
&=-\frac{1}{2}u+\px \left[ xu-2tH\px u+2t|u|^2u \right].
\end{align*}
This motivates the definition
\begin{align*}
L^{NL} u &:= x u - 2t H \partial_x u + 2t |u|^2 u.
\end{align*}
Then,
\begin{align} \label{v LNL}
v &=  \partial_x L^{NL} u-\frac{1}{2} u.
\end{align}
We know $v$ satisfies
\begin{equation} \label{linearized mBO}
    v_t+H\partial_x^2v =  \partial_{x}(2|u|^2v+u^2\bv).
\end{equation}
So plugging in, we get
\begin{align*}
     (\partial_x L^{NL}u-\frac{1}{2}u)_t+H\partial_x^2(\partial_x L^{NL}u-\frac{1}{2}u) &= 2\partial_{x}(|u|^2(\partial_x L^{NL}u-\frac{1}{2}u))+\px(u^2(\px \widebar{L^{NL}u}-\frac{1}{2}\bu)).\\
    -\frac{1}{2}\partial_x (\vert u\vert^2 u)+ \partial_x (\partial_t+H\px^2)L^{NL} u =& -\frac32 \partial_x (\vert u\vert^2 u) + \px(2|u|^2 \px L^{NL}u+u^2\px \widebar{L^{NL}u}) \\
(\partial_t+H\px^2)L^{NL} u =& -  \vert u\vert^2 u + (2|u|^2 \px L^{NL}u+u^2\px \widebar{L^{NL}u}).
\end{align*}
Note that in the last line, when integrating both sides, there is a constant depending on $t$.  However, this constant must be zero, which can be seen by plugging $L^{NL}u$ into the equation.

\subsection{\texorpdfstring{Energy Estimate for $L^{NL}u$}{Energy Estimate for LNLu} }

Now that we have the equation $L^{NL}u$ satisfies, we can compute the energy estimate.  We have
\begin{align*}
    (\partial_t+H\px^2)L^{NL} u =& -  \vert u\vert^2 u + (2|u|^2 \px L^{NL}u+u^2\px \widebar{L^{NL}u}).
\end{align*}
Let $z:=L^{NL}u,$ and assume $t\geq 1.$ Then, 
\begin{align*}
    (\partial_t+H\px^2)z =& -  \vert u\vert^2 u + (2|u|^2 \px z+u^2\px \widebar{z}).
\end{align*}
We calculate the energy estimate in the usual way by multiplying by $\bz,$ taking twice the real part, and integrating.
\begin{align*}
    \frac{d}{dt}\|z\|_{L^2}^2=&2\Rea \int  -  |u|^2 u \bz+ (2|u|^2 \px z+u^2\px \widebar{z})\bz \\
    =&2\Rea \int  -  |u|^2 u \bz+ 2|u|^2 \bz\px z+u^2\bz \px\bz \\
    =&-2\Rea \int |u|^2 u \bz+\Rea \int u^2 \px \bz^2+2 \int |u|^2 \px |z|^2\\
    =&-2\Rea \int |u|^2 u \bz-\Rea \int \px u^2  \bz^2-2 \int\px |u|^2  |z|^2.
\end{align*}
Then, applying the bootstrap assumptions, we have
\begin{align*}
    \left| \frac{d}{dt}\|z\|_{L^2}^2\right|\lesssim& \|u\|_{L^\infty}^2\|u\|_{L^2}\|z\|_{L^2}+\|u\|_{L^\infty}\|u_x\|_{L^\infty}\|z\|_{L^2}^2\\
    \leq& D^2\ep^3 t^{-1}\|z\|_{L^2} + D^2\ep^2t^{-1}\|z\|_{L^2}^2.
\end{align*}
Applying the nonlinear Gr\"onwall inequality, i.e. Lemma \ref{NL Gronwall} with $t_0=1$ we get
\begin{align*}
    \|z\|_{L^2}\lesssim \ep t^{CD^2\ep^2},
\end{align*}
where $C$ is some universal constant.
This proves 
\begin{align*}
     \|L^{NL}u\|_{L^2}\lesssim \ep t^{CD^2\ep^2}.
\end{align*}
Note that 
\begin{align*}
    L^{NL}u=Lu+2t|u|^2u,
\end{align*}
so we can use this to get the corresponding bound for
$\|Lu\|_{L^2}.$  We have
\begin{align*}
    \|Lu\|_{L^2}\lesssim&  \|  L^{NL}u\|_{L^2}+  \|t|u|^2u\|_{L^2}\\
    \lesssim& \|  L^{NL}u\|_{L^2}+ t\|u\|_{L^\infty}^2\|u\|_{L^2}.
\end{align*}

Applying the bootstrap assumptions and the conservation of mass, we have
\begin{align*}
    \|Lu\|_{L^2}
    \lesssim& \ep t^{CD^2\ep^2}+ D^2\ep^3\lesssim \ep t^{CD^2\ep^2},
\end{align*}
as long as $D \ep<1.$

\subsection{Energy estimate for derivative of the vector field}

As was done for DNLS in \cite{byars2024globaldynamicssmalldata}, we also need an energy estimate for $\px Lu$.  We can again do this by first obtaining an estimate for $\px L^{NL}u.$  This is even easier than the previous energy estimate, because we can find the energy estimate for the linearized equation \eqref{linearized mBO} and then use \eqref{v LNL}. First, the estimate for \eqref{linearized mBO} is

\begin{align*}
    \frac{d}{dt}\|v\|_{L^2}^2 =& 2\Rea\int \partial_{x}(2|u|^2v+u^2\bv)\bv\\
    \lesssim& \|u\|_{L^\infty}\|u_x\|_{L^\infty}\|v\|_{L^2}^2\\
    \leq& D^2\ep^2t^{-1}\|v\|_{L^2}^2.
\end{align*}
Applying Gr\"onwall's inequality and using the fact that $v_0=\px(xu_0)-\frac{1}{2}u_0,$ which from \eqref{initials mbo} is bounded $\|v_0\|_{L^2}\leq \ep,$ this becomes
\begin{align*}
    \|v\|_{L^2}\lesssim \ep t^{CD^2\ep^2}.
\end{align*}
Then, \eqref{v LNL} gives
\begin{align*}
        \|\px L^{NL}u\|_{L^2}\leq   \|u\|_{L^2}+    \|v\|_{L^2}\leq \ep \la t \ra^{CD^2\ep^2}.
\end{align*}
Translating back to $Lu,$ we have
\begin{align*}
     \|\px Lu\|_{L^2}\leq&  \|\px L^{NL}u\|_{L^2} + \|t\px (|u|^2u)\|_{L^2}\\
     \leq& \ep \la t \ra^{CD^2\ep^2}+ t\|u_x\|_{L^\infty}\|u\|_{L^\infty}\|u\|_{L^2}\\
     \leq& \ep \la t \ra^{CD^2\ep^2}+ D^2\ep^3\lesssim \ep \la t \ra^{CD^2\ep^2}.
\end{align*}

\section{Closing the bootstrap argument\label{Close bootstrap section}}
Now we close the bootstrap arguments, therefore completing the proof of Lemma \ref{close bootstrap lemma}. We begin by closing the bootstrap argument for $u^\pm.$ 
By combining the difference bound \eqref{spatial diff bnds mbo} and the vector field energy estimate \ref{vector field EE}, we have
\begin{align*}
     \|e^{-i\phi^\pm}u^\pm-t^{-\frac{1}{2}}\gamma^\pm\|_{L^\infty}\lesssim t^{-\frac34}\|Lu\|_{L^2}+|c|t^{-\frac12}\lesssim \epsilon\, \langle t\rangle ^{CD^2\epsilon^2-\frac{3}{4}}+\ep t^{-\frac12}.
\end{align*}
Now, by applying the vector field bound to the $R^\pm$ bound, we have
\begin{align*}
    |R^\pm|\lesssim& \ep t^{-\frac{5}{4}+     CD^2\ep^2}+
   \ep |v|\,t^{-\frac{1}{2}}\la t^{\frac{1}{2}}v\ra^{-N}+ D^2 \ep^3 t^{-\frac54+\delta}\\
    +& D^3 \ep^3 t^{-\frac{3}{2}+\delta}+D^2\ep^3 t^{-\frac{5}{4}+CD^2\ep^2},
\end{align*}
where we showed the second term is integrable in time. 
Then,  \eqref{gamma+ bound} becomes
\begin{align*}
    |\g^{\pm}|\lesssim&  \ep+D^2\ep^3+D^3\ep^3,
\end{align*}
where, in order for $R^{\pm}$ to be integrable in time, we need to ensure that $ CD^2\ep^2 <\frac{1}{4}.$ Now, let $C$ also be the implicit constant from the above equation and the difference bound, so that 
\begin{align*}
    |\g^{\pm}|\leq&C \ep(  1+D^2\ep^2+D^3\ep^2)\\
     \|e^{-i\phi^\pm}u^\pm-t^{-\frac{1}{2}}\gamma^\pm\|_{L^\infty}\leq& C\epsilon\, \langle t\rangle ^{CD^2\epsilon^2-\frac{3}{4}}+C\ep t^{-\frac12}.
\end{align*}
 We need to choose $D$ such that $D>5C$ and $\ep$ small enough such that $CD\ep <\frac{1}{5}.$  Then, using these we have
\begin{align*}
    |\g^{\pm}|\leq&  \ep\lp \frac{D}{5}+ \frac{D}{5}+\frac{D}{10} \rp < \frac{3D}{5}\ep.
\end{align*} 
Then, by triangle inequality,
\begin{align*}
    |u^\pm|\leq& t^{-\frac{1}{2}}|\gamma^\pm| + |e^{-i\phi^\pm}u^\pm-t^{-\frac{1}{2}}\gamma^\pm|\\
   <& \ep t^{-\frac{1}{2}}\left[\frac{3D}{5}\right]+\frac{D}{5}\epsilon\, \la t\rangle ^{CD^2\epsilon^2-\frac{3}{4}}+\frac{D}{5}\ep t^{-\frac12}< D\ep t^{-\frac{1}{2}}.
\end{align*}

The order of choosing $D,\ep$ is important.  First, $C$ is a universal constant which does not depend on $D,\ep,$  and $C>1.$ We choose $D$ large enough so that $D>5C.$  Then choose $\ep <\frac{1}{5CD}$, so that all of the above conditions are met. Note that we will need this stronger condition on $\ep$ to close the second bootstrap argument on $\|u^\pm_x\|_{L^\infty}$. For now, this closes the bootstrap argument for $\|u^\pm\|_{L^\infty}.$

Lastly, we need to close the bootstrap for $\|u_x^\pm\|_{L^\infty}.$ We combine \eqref{spatialu_xbnds+} with the vector field energy estimate \eqref{derivative vector field EE} to get
\begin{align*}
\left\|u^\pm_x(t,vt)\mp \frac{i}{2}t^{-\frac{1}{2}}e^{i\phi^\pm(t,vt)}v\gamma^\pm(t,v)\right\|_{L^\infty}\lesssim& \ep t^{-\frac34+CD^2\ep^2}+
\ep t^{-1}.
\end{align*} 
Now, by applying the vector field bound to the $vR^\pm$ bound, we have
\begin{align*}
    |vR^\pm| \lesssim&D^3\ep^3 t^{-\frac32}+ \ep v^2t^{-\frac{1}{2}}\la t^{\frac{1}{2}}v\ra^{-N}\\
    &+  D^2\ep^3  t^{-\frac54+CD^2\ep^2}+\ep t^{-\frac{7}{4}+CD^2\ep^2}+ D\ep^3 t^{-\frac54+2CD^2\ep^2},
\end{align*}
where we showed the second term is integrable in time. Then,  \eqref{vgamma+ bound} becomes
\begin{align*}
       |v\g^\pm(t)| \lesssim \ep+ D^3 \ep^3 +  D^2\ep^3 + D\ep^3,
\end{align*}
where, in order for $vR^{\pm}$ to be integrable in time, we need $ 2CD^2\ep^2 <\frac{1}{4}.$  This is true since $2CD^2\ep^2 <\frac{2}{25}<\frac{1}{4}.$
 Now, let $C$ also be the implicit constant for the  $v\g^\pm$ bound above, so that 
\begin{align*}
    |v\g^{\pm}|\leq&C \ep(  1 + D^3 \ep^2 +  D^2\ep^2 + D\ep^2).
\end{align*}
Similarly we use $D>5C$, and $CD\ep <\frac{1}{5}.$  Using these we have
\begin{align*}
    |v\g^{\pm}|\leq&  \ep\lp  \frac{D}{5} +  \frac{D}{5} \cdot\frac{1}{25}+ \frac{D}{5} \cdot\frac{1}{25} + \frac{1}{5}\cdot \frac{1}{25}D \rp \\
 <& \frac{2D}{5}\ep.
\end{align*} 
Then, by triangle inequality,
\begin{align*}
    |u^\pm_x|\leq& t^{-\frac{1}{2}}|v\gamma^\pm| + |u^\pm_x(t,vt)\mp \frac{i}{2}t^{-\frac{1}{2}}e^{i\phi^\pm(t,vt)}v\gamma^\pm(t,v)|\\
   <& \ep t^{-\frac{1}{2}}\left[\frac{2D}{5}\right]+\frac{D}{5}\epsilon\, \langle t\rangle ^{CD^2\ep^2-\frac{3}{4}}+\frac{D}{5}\ep t^{-1}< D\ep t^{-\frac{1}{2}},
\end{align*}
This closes the bootstrap argument for $\|u^\pm_x\|_{L^\infty}.$
This completes the proof of Lemma \ref{close bootstrap lemma}, which implies that 
\begin{align*}
    \|u\|_{L^\infty}\leq \|\up\|_{L^\infty}+\|\um\|_{L^\infty}\lesssim \ep t^{-\frac{1}{2}}
\end{align*}
for $t\geq 1.$  We now complete the proof of Theorem \ref{mbo thm} by proving the result in the case of short times.

\begin{proof}[Proof of the case $0\leq t\leq 1$ ]
    Let $D>5$ and $\ep$ small such that $D^2\ep^2<\frac{1}{3}$. We make the bootstrap assumptions:
\begin{align*}
    \|u\|_{L^\infty}\leq D\ep \la t\ra^{-\frac{1}{2}}, \,\,     \|u_x\|_{L^\infty}\leq D\ep \la t\ra^{-\frac{1}{2}}, \,\,\, 0\leq t \leq 1.
\end{align*}
We have the energy estimate
\begin{align*}
    \frac{d}{dt} \|u_x\|_{L^2}^2 =& 2 \int_{\mathbb{R}} \partial_x(|u|^2) |u_x|^2 \, dx + \operatorname{Re} \int_{\mathbb{R}} \partial_x(u^2) \bar{u}_x^2 \, dx\\
    \left| \frac{d}{dt} \|u_x\|_{L^2}^2 \right| \le& 6 \|u\|_{L^\infty} \|u_x\|_{L^\infty} \|u_x\|_{L^2}^2 \le 6 D^2 \epsilon^2 \langle t \rangle^{-1} \|u_x\|_{L^2}^2 \le 6 D^2 \epsilon^2 \|u_x\|_{L^2}^2,
\end{align*}
where we used $\langle t \rangle^{-1} \le 1$ for $t \in [0, 1]$. Applying Gr\"onwall's inequality on $[0, 1]$ yields
\[
\|u_x(t)\|_{L^2}^2 \le \|\px u_{0}\|_{L^2}^2 \exp\left( \int_0^t 6 D^2 \epsilon^2 \, ds \right) \le \epsilon^2 e^{6 D^2 \epsilon^2 t}.
\]
Then we have
\[
\|u_x\|_{L^2} \le  e^{3 D^2 \epsilon^2}  \ep<3\ep. 
\]
By Sobolev embedding, 
\[\|u\|_{L^\infty}\leq \|u\|_{H^1}\leq \|u\|_{L^2}+\|u_x\|_{L^2}<4\ep\lesssim 4\ep \la t\ra^{-\frac{1}{2}}< D\ep \la t\ra^{-\frac{1}{2}}.\] 
This closes the first bootstrap.  For the second, Proposition \ref{vector field estimate general f} gives

\begin{align*}
    \|u_x^{\pm}\|_{L^\infty}\lesssim t^{-\frac{1}{2}}\|u_x^\pm\|_{L^2}^{\frac{1}{2}}\|Lu_x^\pm\|_{L^2}^\frac{1}{2}.
\end{align*}
Lemma \ref{commutator lemma} gives 
$   L(u^\pm_x)=P^\pm L(u_x).$ We also have $    [\px, L]=1,$ so 
\begin{align*}
    \|Lu_x^\pm\|_{L^2}= \|P^\pm Lu_x\|_{L^2}\leq \|Lu_x\|_{L^2}\leq \|\px Lu\|_{L^2}+\|u\|_{L^2}.
\end{align*}
Recall from section \ref{Vector field bnds section} that $v =  \partial_x L^{NL} u-\frac{1}{2} u,$ and 
\begin{align*}
    \frac{d}{dt}\|v\|_{L^2}^2 \leq& D^2\ep^2\la t \ra^{-1}\|v\|_{L^2}^2\leq D^2\ep^2\|v\|_{L^2}^2
\end{align*}
for $0\leq t\leq 1.$  Since $\|v_0\|_{L^2}\leq \epsilon,$ applying Gr\"onwall leads to
\begin{align*}
    \|v\|_{L^2}\leq \ep e^{CD^2\ep^2t} <3 \ep
\end{align*}
for $CD^2\ep^2<1.$ Therefore $\|\partial_x L^{NL} u\|_{L^2}\leq 4\ep.$
Combining all the bounds gives
\begin{align*}
    \|u_x^{\pm}\|_{L^\infty}\lesssim t^{-\frac{1}{2}}(3\ep)^{\frac{1}{2}}(5\ep)^\frac{1}{2}<4\ep t^{-\frac{1}{2}}< D\ep t^{-\frac{1}{2}}.
\end{align*}
This closes the second bootstrap and completes the proof in the range $0\leq t\leq 1.$
\end{proof}
Together, we now have the proof of Theorem \ref{mbo thm}.

\bibliography{refs_combined}{}
\bibliographystyle{plain}
\end{document}